\documentclass[11pt,reqno]{amsart}
\usepackage{color}
\usepackage[active]{srcltx}
\usepackage{a4wide}
\usepackage{amssymb, amsmath, mathtools}
\usepackage{graphicx}
\usepackage{float}
\usepackage[active]{srcltx}
\usepackage[ansinew]{inputenc}
\usepackage{hyperref}

\theoremstyle{plain}
\newtheorem{theorem}{Theorem}[section]
\newtheorem{maintheorem}{Theorem}
\newtheorem{lemma}[theorem]{Lemma}
\newtheorem{proposition}[theorem]{Proposition}

\theoremstyle{remark}

\newtheorem{definition}{Definition}

\newtheorem{remark}[theorem]{Remark}

\numberwithin{equation}{section}

\newcommand{\NN}{{\mathbb{N}}}
\newcommand{\ZZ}{{\mathbb{Z}}}
\newcommand{\RR}{{\mathbb{R}}}

\newcommand{\EU}{{\mathbb{S}}}

\newcommand{\dpt}{\displaystyle}

\usepackage[pagewise]{lineno}

\begin{document}

\title[On the dynamics of rotating rank-one strange attractors families]{On the dynamics of \\ rotating rank-one strange attractors families}

\author[A. A. Rodrigues and B. F. Gon\c{c}alves]{Alexandre A. Rodrigues$^{1}$ and Bruno F. Gon\c{c}alves$^{2}$ \\
 \\ 
$^1$Lisbon School of Economics \& Management, \\Centro de Matem\'atica Aplicada \`a Previs\~{a}o e Decis\~{a}o Econ\'omica \\Rua do Quelhas 6,  Lisboa 1200-781, Portugal \\
$^{2}$Centro de Matem\'atica da Universidade do Porto, \\ \medskip Rua do Campo Alegre s/n, Porto 4169-007, Portugal \\
 }

\thanks{ The first author was supported by the Project CEMAPRE/REM -- UIDB /05069/2020 financed by FCT/MCTES through national funds. The second author was partially supported by CMUP (UID/MAT/00144/2020), which is funded by FCT with national (MCTES) and European structural funds through the programs FEDER, under the partnership agreement PT2020.     }
\email{alexandre.rodrigues@fc.up.pt, arodrigues@iseg.ulisboa.pt,  \quad brunoffg9@gmail.com }


\begin{abstract}
In this article,  we study a two-parameter family of rotating rank-one maps defined on $\EU^1\times [1, 1+b]\times\EU^1$, with $b\gtrsim 0$, whose dynamics  is characterised  by a coupling of a  family of planar maps exhibiting rank-one strange attractors and an Arnold family of circle maps. The main result is about the dynamics on the skew-product, which  is   governed by the existence and prevalence of strange attractors in the corresponding resonance tongues of the Arnold family. 
 The strange attractors carry the unique physical measure of the system, which determines the behaviour of Lebesgue-almost all initial conditions. 
 
  This phenomenon  can be considered as the transition dynamics from a strange attractor with one positive Lyapunov exponent to hyperchaos. 
  Besides an analytical rigorous proof, we illustrate the main results with numerical simulations. We also conjecture how persistent  hyperchaos can be obtained.
 \end{abstract}

\maketitle \setcounter{tocdepth}{1}

\section{Introduction}\label{intro}
The qualitative theory of dynamical systems has seen a  development since the groundbreaking contributions of Poincar\'e and Lyapunov over a century ago. It provided a   framework to describe and understand a wide range of phenomena in several areas as physics, life science and engineering \cite{AHL2001, Chernov93}.
Such a success benefits  from the fact that the law of evolution in various problems is static and does not change with time.   However, many real world problems involve time-dependent parameters and, furthermore, one wants to understand control, modulation or other effects. In doing so, \emph{periodically} or \emph{almost periodically} driven systems are special cases and a general theory for arbitrary time-dependence is desirable. Many of  the  established concepts, methods and results for autonomous systems are not applicable to these cases and require an appropriate extension which fits in the \emph{theory of non-autonomous dynamical systems} \cite{KR11}.

In classical mechanics, dissipative non-autonomous systems received limited attention, in part because it was believed that, for this class of systems, all solutions tend toward Lyapunov stable sets. Evidence that second order equations with a periodic forcing term can have interesting behaviour has appeared in the study of van der Pol's equation, which describes an oscillator with nonlinear damping. Results in Cartwright and Littlewood \cite{CL45} pointed out an attracting set more complicated than a fixed point or an invariant curve. Levinson obtained detailed information for a simplified model \cite{Levinson49}.
Examples from the dissipative category include the equations of Lorenz, Duffing equation and Lorentz gases acted on by external forces \cite{Chernov93}. To date there has been very little systematic investigation of the effects of  time-periodic perturbations, despite being natural for the modelling of many biological effects, see Rabinovich \emph{et al} \cite{Rabinovich06}.

The works by Chen, Oksasoglu and Wang \cite{Chen2013} and Rodrigues \cite{Rodrigues2021} deal with heteroclinic bifurcations in time-periodic perturbations in the dissipative context. The authors have shown, for a set of parameters with positive Lebesgue measure, the existence of an attracting torus, infinitely many horseshoes and strange attractors with one positive Lyapunov exponent.  

In \cite{LR18, TD} the authors analysed a family of periodic perturbations of a weakly attracting robust heteroclinic network defined on the two-sphere. They derived the first return map near the heteroclinic cycle for small amplitude of the perturbing term and, under an open condition in the space of parameters (defined in Theorem 1 of \cite{LR18}), they reduced the analysis of the non-autonomous system to that of a two-dimensional map on a circloid. However, without this assumption, the first return map to a cross section has three-components.  The analysis of a non-autonomous map is, in general, very difficult to tackle.

\subsection*{Strange attractors}
A compact attractor is
said to be \emph{strange} if it contains a dense orbit with at least one positive Lyapunov exponent.
A dynamical phenomenon in a one-parameter family of maps is said to be \emph{persistent} if it occurs for a set of parameters of  positive Lebesgue measure.
Persistence of chaotic dynamics is physically relevant because it means that a given phenomenon is numerically observable with positive probability  \cite{Barrientos_book}.

Strange attractors are of fundamental importance in dynamical systems; they have
been observed and recognized in many scientific disciplines  \cite{AHL2001, Homb2002}.  Atmospheric physics provides one of the most striking   examples of strange attractors observed in natural sciences. We address the reader to \cite{SNN95}  where the authors established the emergence of strange attractors in a low-order atmospheric circulation model.
  Among the theoretical examples that have been studied are the Lorenz and H\'enon attractors, both of which are
closely related to suitable one-dimensional reductions.
The rigorous proof of the strange character of an invariant set is an  involved   challenge and the proof of the persistence (in measure) of such attractors is a challenge.

   In this paper we give a further step towards this analysis.  We provide a criterion for the existence of abundant  strange attractors (in the terminology of \cite{MV93}) near a specific family of diffeomorphisms, using the \emph{Theory of rank-one attractors}  developed by Q. Wang and L.-S. Young \cite{WY2001, WY, WY2003, WY2008}.

\subsection*{Motivation and novelty}
Motivated by the full first return map of  \cite{LR18}  and the Lorenz-84 atmospheric model with seasonal forcing \cite{BSV}, we study a two-parameter family of maps $ \mathcal{F}_{(\varepsilon_1,\varepsilon_2)}$, $\varepsilon_1, \varepsilon_2< 1$, defined on $\EU^1\times [1, 1+b]\times\EU^1$, $b\gtrsim 0$, which can be seen as a coupling of two maps, one exhibiting rank-one strange attractors (in the sense of \cite{WY}) and the other having periodic and quasi-periodic motion (in the sense of Arnold, cf. \cite{MT}). 

The family $ \mathcal{F}_{(\varepsilon_1,\varepsilon_2)}$, $\varepsilon_1, \varepsilon_2< 1$, is related with the Lorenz 84 model described in \cite{BSV}.  In the latter article, the authors consider a non-autonomous perturbed Lorenz model and studied the dynamics associated to a three-dimensional non-autonomous Poincar\'e map which depends on the oscillating part of the forcing.  
 In general, these families may behave periodically, quasi-periodically or chaotically, depending on specific character of the perturbation. 
 
In the present article, we first show that the uncoupled diffeomorphism has a quasi-periodic attractor (Theorem \ref{prop0}). Then we prove that the coupled diffeomorphism $\mathcal{F}_{(\varepsilon_1, \varepsilon_2)}$ has an invariant circle $\mathcal{C}$ of saddle type, such that  its stable and unstable manifolds   are bounded and  such that the orbits of all points within the absorbing domain of   $\mathcal{F}_{(\varepsilon_1, \varepsilon_2)}$ are attracted to the closure of $W^u(\mathcal{C})$ (Theorem \ref{Th2.1}). 
A new class of strange attractors has been found, governed by the dynamics of the unstable manifold of a quasi-periodic orbit which, in some cases, is \emph{ergodic}.   
Using the theory associated to Arnold tongues, we   show how a strange attractor may be strictly contained within the closure of $W^u(\mathcal{C})$ (Theorem \ref{Th2.5}).

 The irreducible strange attractors in the present paper have one direction of instability, are nonuniformly hyperbolic and structurally unstable, although their existence is a prevalent phenomenon.
These results are  expectable  in the transition   from a strange attractor (one positive Lyapunov exponent) to hyperchaos (two positive Lyapunov exponents).

 \subsection*{Structure of the article}
 This article is organised as follows. 
 In Section \ref{s: object} we describe precisely our object of study. In Section \ref{s: main} we state the main results of the article after the introduction of some basic definitions in Section \ref{s: preliminaries}.
 After reviving the Theory of Rank-one strange attractors in Section  \ref{s: theory}, the proof of the main results is performed in Sections \ref{Proof_A}, \ref{s: Th B proof}, \ref{s: Prova Th C} and \ref{s: Prova Th D}.
 We illustrate the main results of the paper with numerics in  Section \ref{s: numerics}.  
  Section \ref{s: discussion} concludes the paper, where we relate the main results of Section   \ref{s: main}  with others in the literature. We also conjecture about the existence of persistent rank-two attractors (under some mild conditions).
  
Throughout this paper, we have endeavoured to make a self contained exposition bringing together all topics related to the proofs. We have stated short lemmas and we have drawn illustrative figures to make the paper easily readable.

\section{Object of study}
\label{s: object}
In this section, we describe the two-parameter map that is the object of study in this research, as well as some terminology and notation. 
 
\subsection{Model under consideration}
In this article, for $b\gtrsim0$ small, we investigate the dynamics of the following map 
 defined on   $M=\EU^1\times [1,1+ b] \times \EU^1$, where $\EU^1 = \RR \pmod{2\pi}$: \\
\begin{equation}
\label{map_general}
\left(\begin{array}{c}x\\ \\ y\\ \\t\end{array}\right)
\mapsto
 \left(\begin{array}{c}
x+ \alpha_1 + \varepsilon_1 \Psi_1(x,y) + \delta_1\ln ( (y-1)+ \varepsilon_1[ \Psi_2 (x,y)+\varepsilon_2\Psi_4 (x,y,t)]) \pmod{2\pi}
\\
\\
1+ ((y-1)+ \varepsilon_1 g(x,y))^\delta
\\
\\
t+\alpha_2+\delta_2   \ln (  \Psi_3(t))\pmod{2\pi}
\end{array}\right).
\end{equation}
under the following assumptions:\\

\begin{description}
 \item[(H1)] $0< \varepsilon_1, \varepsilon_2 < 1$;\\
  \item[(H2)]  $ \alpha_1, \alpha_2 \in [0,2\pi ]$, $\delta_1,\delta_2\in  \RR^+$ and $\delta>1$;\\
\item[(H3)] $\Psi_3: \EU^1 \to \RR^+$ is a $C^1$--map; \\
\item[(H4)]  $\Psi_1, \Psi_2:\EU^1\times [1, 1+b] \to \RR^+$, $\Psi_4:\EU^1\times [1, 1+b]\times \EU^1 \to \RR^+$ are non-constant  $C^3$--maps where $ \ln (  \Psi_2(x,0))\equiv \ln (  \Psi_2(x)) $ has two non-degenerate critical points,  say $c^{(1)}$ and $c^{(2)}$ (\footnote{The main results of this manuscript are valid for a finite number of critical points -- cf. Section 5.3 of \cite{WY}});\\
 \item[(H5)] $\dpt  \sup_{t\in [0,2\pi]} \left|\frac{\Psi_3'(t)}{\Psi_3(t)}\right|<\frac{1}{\delta_2}$;\\
 \item[(H6)]  $g:\EU^1 \times [1, 1+b] \to \RR^+_0$ is $C^3$--map.
  \bigbreak

\end{description}
Define the two-parameter family of   maps $\mathcal{F}_{(\varepsilon_1, \varepsilon_2)}: \EU^1\times  [1,1+ b]  \times \EU^1 \to  \EU^1\times  [1,1+ b] \times \EU^1$ as:\\
$$
\mathcal{F}_{(\varepsilon_1, \varepsilon_2)}(x, y, t)=(\mathcal{F}_1(x,y,t), \mathcal{F}_2(x,y,t), \mathcal{F}_3(x,y,t)),
$$
where:\\
 \begin{equation}
  \left\{ 
\begin{array}{lll}
\mathcal{F}_1(x, y, t)&=&x+ \alpha_1 + \varepsilon_1 \Psi_1(x,y) + \delta_1\ln ( (y-1)+ \varepsilon_1[ \Psi_2 (x,y)+\varepsilon_2\Psi_4 (x,y,t)]) \pmod{2\pi}  \\ \\
\mathcal{F}_2 (x, y, t)&=& 1+((y-1)+ \varepsilon_1 g(x,y))^\delta \\\\
\mathcal{F}_3 (x, y, t)&=& t+\alpha_2+\delta_2   \ln (  \Psi_3(t))\pmod{2\pi}.\\
\end{array}
\right.\\
  \end{equation}
 \bigbreak

The coordinates $(x,y)$ should be understood as polar coordinates ($x$ -- angular component; $y$ --radial component). 
  For $\varepsilon_1 \geq 0$ and $\varepsilon_2=0$, the dynamics of $(\mathcal{F}_1, \mathcal{F}_2)$ and $\mathcal{F}_3$ are uncoupled (independent).    For the sake of simplicity, we denote the map $(\mathcal{F}_1, \mathcal{F}_2)$ by  $\mathcal{T}_{(\varepsilon_1, \varepsilon_2)}$.\\
 
\begin{remark}
\label{cont1}
Since $y\in  [1,1+ b] $ and $\delta>1$ then,  for small $\varepsilon_1,\varepsilon_2 \geq 0$,  the map $\mathcal{T}_{(\varepsilon_1, \varepsilon_2)}$ 
is dissipative near the circle $y=1$ (\footnote{``Dissipative'' in the sense that the product of the eigenvalues of $D\mathcal{T}_{(\varepsilon_1, \varepsilon_2)}$, at all points of its domain,  has modulus less than 1.}). Indeed, for the radial component we may write:
$$
\frac{\partial \mathcal{F}_2}{\partial y}(x,y)= \delta [(y-1)+\varepsilon_1g(x,y)]^{\delta-1}.\left[1+\varepsilon_1 \frac{\partial g}{\partial y}(x,y)\right] =\mathcal{O}(\varepsilon_1),
$$
where $\mathcal{O}$ stands for the standard Landau notation.  \\
\end{remark}

Let us call by $\mathcal{A}\subset  \EU^1\times  [1,1+ b] $  the \emph{absorbing domain} of $\mathcal{T}_{(\varepsilon_1, \varepsilon_2)}$ -- orbits starting in $\mathcal{A}$ stay there for all positive iterates. \\

\begin{remark}
\label{injective}
The positive map $\mathcal{F}_3$ only depends on $t$ and   may be seen as an injective map on the circle. Indeed, by \textbf{(H5)}, one gets:
        \begin{eqnarray*}
\sup_{t\in [0,2\pi]} \left|\frac{\Psi_3'(t)}{\Psi_3(t)}\right|<\frac{1}  {\delta_2}&\Rightarrow & \forall t\in [0,2\pi], \quad  -\frac{1}{\delta_2}< \frac{\Psi_3'(t)}{\Psi_3(t)} < \frac{1}{\delta_2} \\
& \Leftrightarrow& \forall t\in [0,2\pi], \quad 0< 1+\delta_2\frac{\Psi_3'(t)}{\Psi_3(t)} < 2\\
& \Rightarrow& \forall t\in [0,2\pi],  \quad \frac{\partial \mathcal{F}_3(x,y,t)}{\partial t} >0\\
& \Leftrightarrow& \mathcal{F}_3(x,y,t) \text{ is injective in } t.\\
      \end{eqnarray*}
      
 In particular, we    are not allowing  the existence of homoclinic orbits for  $ \mathcal{F}_3$ \cite{Aronson}. \\
\end{remark}

 \subsection{Remarks}
 \label{imp_remarks}
 With respect to the map \eqref{map_general}, we would like to point out some remarks.\\
  \begin{itemize}

 \item Hypotheses \textbf{(H1)--(H2)} are technical on the theory of normally hyperbolic perturbations and the condition $\delta>1$ guarantees the dissipativeness of   $ \mathcal{T}_{(\varepsilon_1,\varepsilon_2)}$; Hypotheses \textbf{(H3)--(H6)} will be necessary to apply the rank-one strange attractors theory \cite{WY}. \\

 \item The term $((y-1)+ \varepsilon_1 g(x,y))^\delta$ in $\mathcal{F}_2 $ plays the role of the ``contracting'' radial  term necessary to compute the limit family of $ \mathcal{T}_{(\varepsilon_1,0)}$.\\
 \item The family $ \mathcal{T}_{(\varepsilon_1,0)}$, $\varepsilon_1 \gtrsim 0$, may be seen as the first return map near a Bykov attractor to a cross section transverse to a cycle where  the two-dimensional manifolds unfold  from the coincidence (cf. \cite{Rodrigues2022}).\\
  \end{itemize}
 
 Since the family of maps $ \mathcal{T}_{(\varepsilon_1,0)}$, $\varepsilon_1 \gtrsim 0$, exhibits rank-one strange attractors in a \emph{persistent} way  and  $\mathcal{F}_3$ is conjugate to a rigid rotation map, we decided to call the resulting dynamics as ``\emph{rotating rank-one strange attractors families}''. This justifies the title of the present manuscript; the theoretical analysis of the map \eqref{map_general} is the goal of the article.
 
  \section{Preliminaries}
  \label{s: preliminaries}
   In this section, we introduce useful terminology that will be used in the rest of the paper.
   
   We assume that $M$ is the set  defined by $\EU^1\times [1,1+b]\times \EU^1$, $b\gtrsim 0$, endowed with the $C^r$--usual norm $\|\star\|_{\mathbb{C}^r}$ in the quotient space, $r >1$.

If $A\subset M$, let  ${\rm int}(A)$ and $\overline{A}$ be the topological interior and closure of $A$, respectively.  We denote by $\ell_2$ and ${\rm dist}_2$ the  induced Lebesgue measure and distance in $\EU^1\times [1,1+b]$,  respectively.  Analogously, we define  $\ell_3$ and ${\rm dist}_3$ as the  induced Lebesgue measure and distance in $M$.

     \begin{definition}
  Let $F : M \to M$ be a $C^3$--diffeomorphism, $z=(x,y,t)\in M$ and $A\subset M$ such that $F(A)\subset A$.\\
  
    \begin{enumerate}
  \item The \emph{forward orbit} of $z$ under $F$ is the set $\{F^j(z), j\in \NN_0\}$ and is denoted by ${\rm Orb}_F(z)$. \\ 
  
  \item The \emph{omega-limit set} of $z \in M$ is the set $\dpt \bigcap_{n\in \NN}\overline{\{F^k(z), k\geq n\}}$ and is denoted by $\omega_F(z)$;\\
  
  \item The \emph{basin of attraction} of $A$ for $F$ is the set $\{z \in M: \omega_F(z)\subset A\}$ and it is denoted by $\mathcal{B}_F(A)$. \\
  
    \item The $F$--invariant set $A \subset M$ is called \emph{topologically transitive} if there exists a point $z \in A$ such that   ${\rm Orb}_F(z)$ is dense in $A$ ($\Leftrightarrow \overline{{\rm Orb}_F(z)}=A$); \\

  \item The $F$--invariant compact set $A \subset M$    is called a \emph{strange attractor} if there exists $z\in M$ such that: \\
  \begin{enumerate}
  \item  $\overline{{\rm Orb}_F(z)}=A$;\\
  \item there exists  $0\neq v \in T_z M$  where $\|DF^n(z)v\|\geq k\lambda ^n$, for all $n\in \ZZ$, and some $k>0$ and $\lambda>1$;\\
   \item  the  basin of attraction of $A$ has nonempty interior (in particular, it contains a set with positive Lebesgue measure).\\
    \end{enumerate}

  \end{enumerate}
 
     \end{definition}
 
  As usual, similar definitions can be given with $F$ replaced by the power $F^m$, with $m\in \NN$.  When there is no risk of misunderstanding, we omit the subscript $F$ in ${\rm Orb}_F$, $\omega_F$ and $\mathcal{B}_F$.
    
 \begin{definition}
We say that $ F $ possesses a \emph{strange attractor supporting an ergodic SRB measure} $\nu$ if: \\
\begin{enumerate}
\item $ F $ exhibits a strange attractor $\Omega$ in a given $F$--invariant open set $U\subset M$:
$$
\Omega = \bigcap_{m=0}^{+\infty} F^m(\overline{U});
$$
 
\item the conditional measures of $\nu$ on unstable manifolds are equivalent to the Riemannian volume on these leaves and \\
\item  for Lebesgue almost points $z\in U\subset M$ and  for every continuous function $\varphi : U\rightarrow \RR$, we have:
\begin{equation}
\label{limit2}
\lim_{n\in \NN} \quad \frac{1}{n} \sum_{i=0}^{n-1} \varphi \circ F ^i(z) = \int  \varphi \, d\nu.
\end{equation}

\end{enumerate}
 
\end{definition}
   
  The measure $\nu$ is known as a Sinai-Ruelle-Bowen measure (SRB measure). It is natural to link sets of positive Riemannian volume with observable events. If we do so, then the SRB measure $\nu$ is observable because temporal and spatial averages coincide for a set of initial data of full Riemannian volume in the basin.

\section{Main results}
\label{s: main}
Our first result is about the dynamics of $\mathcal{T}_{(\varepsilon_1,0)}$ with $\varepsilon_1 \gtrsim 0$. Theorem \ref{th:transv} provides the existence of a strange attractor for a set of parameters with positive Lebesgue measure, and a nontrivial basin of attraction for  it.
\begin{maintheorem} 
 \label{th:transv}
For $\varepsilon_2=0$, then there exists $\delta_1^\star> 0$ and $ 0<\varepsilon_1^\star\ll 1$ such that if $\delta_1\gg\delta_1^\star$,   there exists  $\mathcal{G}_1\subset [0, \varepsilon_1^\star]$ with positive Lebesgue measure such that for every $\varepsilon_1\in\mathcal{G}_1$ the map $\mathcal{T}_{(\varepsilon_1,0)}$  exhibits an irreducible strange attractor $\Omega$ that supports a unique ergodic SRB measure $\nu$.  
 
   \end{maintheorem}
   The strange attractor $\Omega$ shadows the entire circle defined by $y=1$.
    The proof of this result relies on the reduction of  $\mathcal{T}_{(\varepsilon_1,0)}$  to a Misiurewicz map and is performed in Section   \ref{Proof_A}, using the Theory of Rank-one attractors developed by Wang and Young \cite{WO, WY2001, WY2003, WY2008}.  Furthermore the theory allows to conclude that:  
 
    \begin{enumerate}
    \item   if $\varepsilon_1\in [0, \varepsilon_1^\star] $,  then the map   $\mathcal{T}_{(\varepsilon_1,0)}$ has a hyperbolic periodic point $p_a\in \EU^1\times [1,1+b]$ (of saddle-type) having a transverse homoclinic connection.  The transversality of the invariant manifolds is a property that holds in open sets of the parameter space. 
    \item  if $\varepsilon_1\in \mathcal{G}_1$, then $\overline{W^u(p_a)}$ does not escape from $\mathcal{A}$ and  generates the irreducible strange attractor $\Omega$. Attracting periodic attractors cannot be present in the dynamics.\\
    \end{enumerate}

    \begin{remark}
       \label{mean_a}
   The set $\mathcal{G}_1$ depends on $\delta_1$.     The subscript $a\in [0,2\pi)$ of the previous Remark relies on the Misiurewicz map and depends on $\varepsilon_1>0$ and $\delta_1$ in the following way (cf Section \ref{Proof_A}):
    $$\varepsilon_1 = \exp\left(-\frac{a+2n\pi}{\delta_1}\right)>0, n\in \NN.$$
\end{remark}
The strange attractor $\Omega$ stated in Theorem \ref{th:transv} is non-uniformly hyperbolic because it contains \emph{critical points} -- points belonging to a dense orbit for which a nonzero tangent vector exists which is contracted both by positive and negative iterates of $\mathcal{T}_{(\varepsilon_1,0)}$.
 With respect to the map $\mathcal{T}_{(\varepsilon_1,0)}$, we are going to concentrate the study of the dynamics on $\delta_1\gg1$ and $\varepsilon_1 \in \mathcal{G}_1$, where we can observe  abundance of irreducible rank-one strange attractors. \\


  \begin{figure}[h]
\begin{center}
\includegraphics[height=7cm]{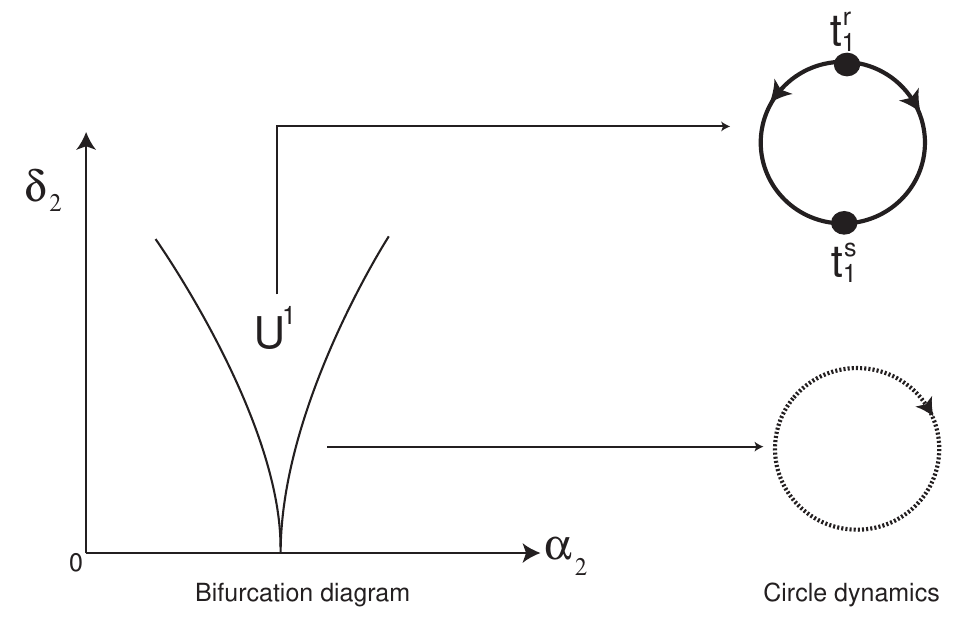}
\end{center}
\caption{\small    Arnold's tongue with rotation number 1 in the $(\alpha_2, \delta_2)$--parameter plane.  The closure of the unstable manifold of the saddle fixed point $t_1^r$ is $\EU^1$. Outside resonance wedges, the circle $\EU^1$ is the minimal attractor. }
\label{Bruno1.1}
\end{figure}
 
  Let $p, q\in \NN$ such that $p\leq q$ and $p,q$ are  relatively prime.
Let $\mathcal{U}^{p/q}$ be the  Arnold tongue associated to the map $\mathcal{F}_3$ with rational rotation number $p/q$ (cf.  \cite{Aronson, MT, SST}). As illustrated in Figure \ref{Bruno1.1}, for $\delta_2>0$ satisfying \textbf{(H5)}, the dynamics of $\mathcal{F}_3$ is characterized by:
\begin{itemize}
\item if $(\alpha_2, \delta_2) \in {\rm int }\,  (\mathcal{U}^{p/q})$, then the dynamics is governed by an even (finite) number of sinks and sources;
\item if $(\alpha_2, \delta_2) \notin \mathcal{U}^{p/q}$, then the dynamics is conjugated to an irrational rotation.  The entire circle $\EU^1$ is the minimal set in the sense that it cannot be decomposable. \\
\end{itemize}

 For $\varepsilon_1 \gtrsim 0$ fixed and $\varepsilon_2=0$,  the hyperbolic periodic orbit $p_a$ of $\mathcal{T}_{(\varepsilon_1,0)}$, guaranteed by Theorem \ref{th:transv}, corresponds to a normally hyperbolic invariant circle of $\mathcal{F}_{(\varepsilon_1,\varepsilon_2)}$ defined by $\mathcal{C}_0=\{p_a\}\times \EU^1$, which persist under $C^3$--perturbations, as well as its stable and unstable manifolds. Set  $\mathcal{A}'=\mathcal{A}\times \EU^1$ as the absorbing domain of $\mathcal{F}_{(\varepsilon_1,0)}$.\\
 
The next result ensures the existence of \emph{quasi-periodic strange attractors} for $\mathcal{F}_{(\varepsilon_1,0)}$    in the sense that they exhibit a strange attractor with one positive and one zero Lyapunov expoent. The proof is performed in Section \ref{s: Th B proof}.

\begin{maintheorem}  
\label{prop0}
If   $0<\varepsilon_1< \varepsilon_1^\star$,  $\delta_1 \gg1$, $\varepsilon_1\in \mathcal{G}_1$ and $(\alpha_2, \delta_2) \notin \mathcal{U}^{p/q}$,  then there exists a $\mathcal{F}_{(\varepsilon_1,0)}$--invariant circle $\mathcal{C}_0=\{p_a\}\times \EU^1$  of saddle-type  such that:\\
\begin{enumerate}
 \item there exists a point $x\in \mathcal{A}'$ whose $\mathcal{F}_{(\varepsilon_1, 0)}$--orbit is dense in   $\overline{W^u(\mathcal{C}_0)}$.\\
\item there exist  $z\in \mathcal{A}'$, constants $k>0, \lambda>1$ and a vector $v\in T_zM$ such that:\\
  $$\left\|D\mathcal{F}_{(\varepsilon_1,0)}^n(z)v\right\|\geq k\lambda ^n,\quad \forall n\in \NN_0.$$
 \end{enumerate}
 \end{maintheorem}  
 \bigbreak
By construction, for $\varepsilon_1\in \mathcal{G}_1$, we have  $\overline{W^u(\mathcal{C}_0)} =\overline{W^u(p_a)}\times \EU^1$.
 For $\varepsilon_1>0$ and  $\varepsilon_2 \gtrsim 0$, we denote by  $\mathcal{C}_{\varepsilon_2}$ the hyperbolic continuation of  $\mathcal{C}_0$.  
    The next result says that there exists chaos in the dynamics associated to  \eqref{map_general}; the circle $\mathcal{C}_{\varepsilon_2}$ is quasi-periodic ($\Rightarrow$ topologically transitive) for a set of parameter values having positive   Lebesgue measure. If $A,B$ are two submanifolds of $M$, the notation $A\pitchfork B$ means that $A$ and $B$ meets transversally. \\

 \begin{maintheorem}
    \label{Th2.1}
  With respect to $\mathcal{F}_{(\varepsilon_1,\varepsilon_2)}$, if   $\delta_1 \gg1$ and  $0<\varepsilon_1< \varepsilon_1^\star$,   then there exists $\varepsilon^\star_2>0$ such that:\\
  \begin{enumerate}
 \item    for all $\varepsilon _2\in [0, \varepsilon^\star_2]$, there exists a non-empty open set $U' \subset M $ such that if $(x,y,t)\in U'$, then $
\omega(x, y, t) \subset\overline{W^u(\mathcal{C}_{ \varepsilon_2})};$ 

  \item   for all $\varepsilon _2\in [0, \varepsilon^\star_2]$ there exists a hyperbolic $C^3$--horseshoe $\mathcal{H}_{\varepsilon_2}$ dynamically defined as $$\overline{W^u(\mathcal{C}_{ \varepsilon_2})\pitchfork {W^s(\mathcal{C}_{ \varepsilon_2}})};$$

\item   there exists a set $\Gamma_{\varepsilon_2}\subset [0, \varepsilon^\star_2]$ with positive Lebesgue measure  such that if $\varepsilon_2\in \Gamma_{\varepsilon_2}$ then the restriction of $\mathcal{F}_{(\varepsilon_1, \varepsilon_2)}$ to the circle $\mathcal{C}_{\varepsilon_2}$ is smoothly conjugate to an irrational rotation. \\
\end{enumerate}
 
\end{maintheorem}  
 The proof of this result is done in Section \ref{s: Prova Th C}.
  Theorem \ref{Th2.1} does not provide a ``minimal'' description of the attracting set associated to  $\dpt \bigcap_{m=0}^{+\infty}  \mathcal{F}_{(\varepsilon_1, \varepsilon_2)}^m(\mathcal{A}')$ by the following reasons:
    \begin{itemize}
  \item for $\varepsilon _2\in [0, \varepsilon^\star_2]$, the set $\mathcal{C}_{\varepsilon_2}$ may be decomposable;
  \item the hyperbolic $C^3$--horseshoe $\mathcal{H}_{\varepsilon_2}$ is invisible in terms of Lebesgue measure \cite{Bowen};
  \item although the set $\overline{W^u(\mathcal{C}_{ \varepsilon_2})}$  attracts an open set of initial conditions, the associated non-wandering set might be non-topologically transitive  \cite{OW2010}.
    \end{itemize}
    \bigbreak

                In what follows,   $(\alpha_2, \delta_2) \in {\rm int} \, ({\mathcal{U}^{p/q}})$ is such that the dynamics of the Arnold map $\mathcal{F}_3$ is Morse-Smale: there exist periodic points $t^s_1, ..., t^s_q$ and $t^r_1, ..., t^r_q$ in $\EU^1$, such that $t^s_i$ is attracting and $t^r_i$  repelling,   $i\in \{1, ..., q\}$. See an illustration of $\mathcal{U}^1$ in  Figure \ref{Bruno1.1}. \\
                
The next result shows that, under precise conditions on $(\alpha_2, \delta_2) $,   the set of parameters $(\varepsilon_1, \varepsilon_2)$ for which $\mathcal{F}_{(\varepsilon_1, \varepsilon_2)}$ has a strange attractor, has positive Lebesgue measure. 

  \begin{maintheorem}
   \label{Th2.5}
   If $\varepsilon_1<\varepsilon_1^\star$, then for all  $\dpt 
 (\alpha_2, \delta_2) \in {\rm int} \, ({\mathcal{U}^{p/q}})$ there exists a set $\tilde{\mathcal{G}}_1\subset [0, \varepsilon_1^\star]$ with positive Lebesgue measure and $ \varepsilon_2^\star>0$ such that
  if $\varepsilon_2< \varepsilon_2^\star$ and $ \varepsilon_1\in \tilde{\mathcal{G}}_1 $, then the map $\mathcal{F}_{(\varepsilon_1, \varepsilon_2)}$ has an irreducible strange attractor contained in  $\overline{W^u( p_a^\star)}$, where $p_a^\star$ is a periodic point of saddle-type.
 
 \end{maintheorem}  
 The strange attractor of Theorem \ref{Th2.5} shadows the set $\EU^1\times \{1\}\times {\rm Orb}_{\mathcal{F}_3}(t_1^s)$.
The proof of this result is the main novelty of this article and is performed in Section \ref{s: Prova Th D}. In general, the set $\mathcal{G}_1$ of Theorem \ref{th:transv} is different from $\tilde{\mathcal{G}}_1$ of Theorem \ref{Th2.5}. The strange attractors stated in Theorem \ref{Th2.5} contains, but do not coincide with the hyperbolic horseshoes of  Theorem \ref{Th2.1}.

\section{Theory of rank-one strange attractors: \\ a brief overview}
\label{s: theory}

We gather in this section a collection of technical facts used later. In what follows, let us denote by $C^2(\EU^1,\RR) $ the set of $C^2$--maps from $\EU^1$ (unit circle) to $\RR$. For $h\in  C^2(\EU^1,\RR) $, let
$$  C(h)= \{x \in \EU^1 : h'(x) = 0\}$$
be the \emph{critical set} of $h$. For $\delta>0$, let $C_\delta$ be the $\delta$--neighbourhood of $C(h)$ in $\EU^1$ and let $C_\delta(c)$ be the $\delta$--neighbourhood of $c\in C(h)$. The terminology $\rm{dist}_1$ denotes the euclidian  metric on $\RR$.

\subsection{Misiurewicz-type maps}
\label{Misiurewicz-type map}
Following Ott and Wang \cite{OW2010}, we say that $h\in  C^2(\EU^1,\RR) $ is a \emph{Misiurewicz-type map}   if the following assertions hold:
\bigbreak
\begin{enumerate}
\item There exists $\delta_0>0$ such that: \\
\begin{enumerate}
\item $\forall x \in C_{\delta_0}$, $h''(x)\neq 0$ and\\
\item $\forall c\in C(h)$ and $n\in \ZZ^+$, ${\rm{dist}_1}(h^n(c), C(h))\geq \delta_0$.\\
\end{enumerate}

\bigbreak

\item There exist constants $b_0, \lambda_0 \in \RR^+$ such that for all $\delta<\delta_0$ and $n\in \NN$, we have: \\
\begin{enumerate}
\item if $h^k(x)\notin C_\delta$ for $k\in\{0, ,..., n-1\}$, then $|(h^n)'(x)|\geq b_0\,  \delta\,  \exp(\lambda_0\, n)$.\\
\item  if $h^k(x)\notin C_\delta$ for $k\in\{0, ,..., n-1\}$ and $h^n(x)\in C_{\delta_0} $, then $|(h^n)'(x)|\geq b_0\,     \exp(\lambda_0\, n)$.\\
\end{enumerate}
\bigbreak
\end{enumerate}

\begin{figure}[h]
\begin{center}
\includegraphics[height=7cm]{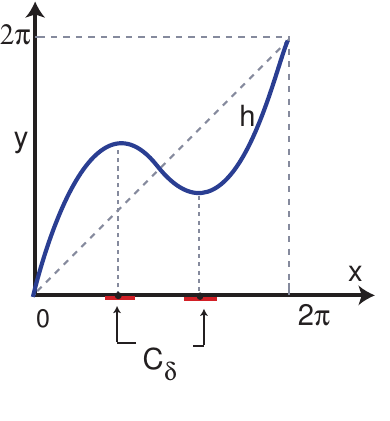}
\end{center}
\caption{\small  ``Shape'' of a Misiurewicz-type map $h:\EU^1 \rightarrow \EU^1$ with two critical points. For $\delta>0$, the set  $C_{\delta}$  is a neighbourhood of the set of critical points $C(h)$.    }
\label{misiurewicz1}
\end{figure}

For $\delta>0$, the set  $\EU^1$  may be divided into two regions: $C_{\delta}$ and $\EU^1 \backslash C_{\delta}$. In $\EU^1 \backslash C_{\delta}$, $h$ is essentially uniformly expanding; in $C_{\delta}\backslash C$,  although $|h '(x)|$ is small, the orbit of $x$ does not return to $C_{\delta}$   until its derivative has regained an amount of exponential growth. We suggest the reader to observe Figure \ref{misiurewicz1} to have in mind the shape of a Misiurewicz-type map.

\subsection{``Rank-one'' maps}
\label{rank_one}

Let $N= \EU^1 \times [0,1]$, induced with the usual topology. We consider the two-parameter family of maps $F_{(a,b)}: N \rightarrow N$, where $a \in[0, 2\pi]$ where $0$ and $2\pi$ are identified, and $b \in \RR$ is a scalar (\footnote{Although related, the two-parameter map $F_{(a,b)}$ is not the same as   $\mathcal{F}_{(\varepsilon_1, \varepsilon_2)}$ of \eqref{map_general}.}). Let $B_0 \subset \RR\backslash\{0\}$ with  0 as an accumulation point. Rank-one theory asks   the following hypotheses:
\bigbreak

\begin{description}
\item[\text{(P1) Regularity conditions}]
\begin{enumerate}
\item For each $b\in B_0$, the function $$(x,y,a)\mapsto F_{(a,b)}(x,y)$$ is at least $C^3$--smooth.
\item Each map $F_{(a,b)}$ is an embedding of $N$ into itself (ie, $F_{(a,b)}(N)\subset N$).
\item There exists $k\in \RR^+$ independent of $a$ and $b$ such that for all $a \in  [0,2\pi]$, $b\in B_0$ and $(x_1, y_1), (x_2, y_2) \in M$, we have:
$$
\frac{|\det DF_{(a,b)}(x_1, y_1)|}{|\det DF_{(a,b)}(x_2, y_2)|} \leq k.
$$
\end{enumerate}

\bigbreak
\item[\text{(P2) Existence of a singular limit}]
For $a\in [0,2\pi]$, there exists a map $$F_{(a,0)}: N \rightarrow \EU^1 \times \{0\}$$ such that the following property holds: for every $(x, y) \in N$ and $a\in [0,2\pi]$, we have
$$
\lim_{b \rightarrow 0} F_{(a,b)}(x,y) = F_{(a,0)}(x,y).
$$

\bigbreak
\item[\text{ (P3) $C^3$--convergence to the singular limit}]   For every choice of $a\in [0,2\pi]$, the maps $(x,y,a) \mapsto F_{(a,b)}$  converge in the $C^3$--topology to $(x,y,a) \mapsto F_{(a,0)}$ on $N\times [0,2\pi]$, as $b$ goes to zero.

\bigbreak
\item[\text{(P4) {\small Existence of a sufficiently expanding map within} the singular limit}]

The\-re exists $a^\star \in [0, 2\pi]$  such that $h_{a^\star}(x)\equiv F_{(a^\star, 0)}(x,0)$ is a Misiurewicz-type map in the sense of Subsection \ref{Misiurewicz-type map}.

\bigbreak
\item[\text{(P5) Parameter transversality}] Let $  C(h_{a^\star})$ denote the critical set of a Misiu\-rewicz-type map $h_{a^\star}$.
For each $x\in C(h_{a^\star})$, let $p = h_{a^\star}(x)$, and let $ {x(\tilde{a})}$ and $ {p(\tilde{a})}$ denote the
continuations of $x$ and $p$, respectively, as the parameter $a$ varies around $a^\star$. The point $ {p(\tilde{a})}$  is the unique point such that $ {p(\tilde{a})}$  and $p$ have identical symbolic itineraries under $h_{a^\star}$ and $h_{\tilde{a}}$, respectively. We have:
$$
\frac{d}{da} h_{\tilde{a}}(x(\widetilde{a}))|_{a=a^\star} \neq \frac{d}{da} p(\tilde{a})|_{a=a^\star}.
$$

\bigbreak

\item[\text{(P6) Nondegeneracy at turns}] For each $x\in C(h_{a^\star})$, we have
$$
\frac{d}{dy} F_{(a^\star,0)}(x,y) |_{y=0} \neq \overline{0}.
$$

\bigbreak

\item[\text{(P7) Conditions for mixing}] If $J_1, \ldots, J_r$ are the intervals of monotonicity of $h_{a^\star}$, then:
\medbreak
\begin{enumerate}
\item $\exp(\lambda_0/3)>2$ (see the meaning of $\lambda_0$ before) and
\medbreak
\item if $Q=(q_{im})$ is the matrix of all possible transitions  defined by:
\begin{equation*}
\left\{
\begin{array}{l}
1 \qquad \text{if} \qquad J_m\subset h_{a^\star} (J_i)\\
0 \qquad \text{otherwise},\\
\end{array}
\right.
\end{equation*}
then there exists $p\in \NN$ such that $Q^p>0$ (\emph{i.e.} all entries of the matrix $Q^p$, endowed with the usual product,  are positive).
\bigbreak
\end{enumerate}

\end{description}

\begin{definition}
\label{terminologia1}
Identifying $\EU^1 \times \{0\}$ with $\EU^1$, we refer to $F_{(a,0)}$  as the restriction $h_a : \EU^1 \rightarrow \EU^1$ defined by $h_a(x) = F_{(a,0)}(x,0)$. This is the \emph{singular limit} of $F_{(a,b)}$.\\
\end{definition}

\subsection{ Rank-one strange attractors}

For attractors with strong dissipation and one direction of instability, Wang and Young conditions \textbf{(P1)--(P7)} are relatively simple and checkable; when satisfied, they guarantee the existence of  strange attractors with a ``package'' of statistical and geometric properties as follows:

\begin{theorem}[\cite{WY}, adapted]
\label{th_review}
Suppose the family $F_{(a,b)}$ satisfies \textbf{(P1)--(P7)}. Then, for all sufficiently small $b\in  B_0$, there exists a subset $\Delta \subset  [0, 2\pi]$ with positive Lebesgue measure such that for $a\in \Delta$, the map $F_{({a},b)}$ admits an irreducible strange attractor $\tilde{\Omega}\subset \Omega$  that supports a unique ergodic SRB measure $\nu$. The orbit of Lebesgue almost all points in $\tilde{\Omega} $ has a positive Lyapunov exponent  and is asymptotically distributed according to $\nu$.
\end{theorem}

The strange attractor $\Omega$ has one direction of instability; this is why the terminology \emph{``rank-one strange attractor''}.
The theory described in \cite{WY2001, WY, WY2003} is
general, in the sense that the conditions under which it is valid depend only on certain
general properties of the   maps and not on specific formulas or contexts.
 \medbreak
 \begin{proposition}[\cite{WY, WY2003}, adapted]
\label{Prop2.1WY}
 Let $A: \EU^1 \rightarrow \RR$ be a $C^3$--map with nondegenerate critical points. Then there exist $L_1$ and $\delta$ depending on $A$ such that if $L \geq  L_1$  and $B: \EU^1 \rightarrow \RR$ is a $C^3$
map with $\|B\|_{\textbf{C}^2} \leq \delta$ and $\|B\|_{\textbf{C}^3}\leq 1$,
 then the family of maps
 $$ t \mapsto t+a+L(A(t) +B(t)), \qquad a\in [0,2\pi[, \qquad t\in \EU^1$$
satisfies \textbf{(P4)} and \textbf{(P5)}. If $L$ is sufficiently large, then  Hypothesis \textbf{(P7)} is also verified.
\end{proposition}


  \section{Proof of Theorem \ref{th:transv}}
  \label{Proof_A}

\subsection{Limit family of $\mathcal{T}_{(\varepsilon_1,0)}$}
\label{s: limit family}
First of all, notice that the domain of  $\mathcal{T}_{(\varepsilon_1,0)}$ is a circloid, the appropriate set where the theory of \cite{WY} may be applied. 
 \subsubsection{Change of coordinates}
\label{change_of_coordinates}
For $\varepsilon_1 \in\,\,  ]0, \varepsilon_1^\star[ $ (small) fixed and $(x,y) \in \mathcal{A}$, let us make the following change of coordinates:
  \begin{equation}
  \label{change1}
  \overline{x} \mapsto {x} \qquad \text{and} \qquad \overline{y} \mapsto \frac{y-1}{\varepsilon_1}.
  \end{equation}
Taking into account that:\\
\begin{eqnarray*}
\mathcal{F}_1 (x,y) &=& x+\alpha_1+{\varepsilon_1} \Psi_1(x,y) + \delta_1  \ln (y-1+{\varepsilon_1}\Psi_2(x,y)) \pmod{2\pi}\\
&=&   x+\alpha_1+{\varepsilon_1} \Psi_1(x,y) + \delta_1  \ln \left[{\varepsilon_1} \left(\frac{y-1}{\varepsilon_1 }+\Psi_2(x,y)\right)\right] \pmod{2\pi}\\
 &=&   x+\alpha_1+{\varepsilon_1}  \Psi_1(x,y) + \delta_1  \ln \varepsilon_1 +\delta_1 \ln  \left[ \frac{y-1}{\varepsilon_1}+\Psi_2(x,y)\right] \pmod{2\pi}
  \end{eqnarray*}
  and
    \begin{eqnarray*}
\mathcal{F}_2 (x,y) &=&1+ ((y-1)+\varepsilon_1g(x,y))^\delta =  1+\varepsilon_1^\delta \left(\frac{y-1}{\varepsilon_1}+ g(x,y)\right)^\delta,
  \end{eqnarray*}
  we may write (in the new coordinates):\\
\begin{eqnarray*}
\mathcal{F}_1  (x,\overline{y}) &=& x+\alpha_1+{\varepsilon_1}  \Psi_1(x,\overline{y}) + \delta_1  \ln \varepsilon_1 +\delta_1 \ln  \left[\left( \overline{y}+\Psi_2(x,\overline{y})\right)\right] \pmod{2\pi}\\ \\
\mathcal{F}_2 (x,\overline{y}) &=&  \varepsilon_1^{\delta-1} \left( \overline{y}+ g(x,\overline{y})\right)^\delta. \\
  \end{eqnarray*}

\subsubsection{Reduction to a singular limit}
\label{ss: reduction}
In this subsection, we compute the singular limit of $\mathcal{T}_{(\varepsilon_1,0)}=(\mathcal{F}_1,\mathcal{F}_2)$ written in the coordinates $(x,\overline{y})$ studied in Subsection \ref{change_of_coordinates}, for $\varepsilon_1 \in\,\, ]0, \varepsilon_1^\star[$.
Let  $k: \RR^+ \rightarrow \RR$ be the invertible map defined by $$k(x)= \delta_1 \ln (x).$$
   Define now the decreasing sequence $(\varepsilon_{n})_n$ such that, for all $n\in \NN$, we have:\\
\begin{enumerate}
\item  $\varepsilon_n\in\, ]0, \varepsilon_1^\star[$ and \\
\item $k(\varepsilon_n) \equiv 0 \pmod{2\pi}$.\\
\end{enumerate}
\medbreak
 
Since $k$ is an invertible map,  for $a \in \EU^1 \equiv [0, 2\pi[$ fixed and $n\geq n_0\in \NN$, let
\begin{equation}
\label{sequence1}
\varepsilon_{( a,n)}= k^{-1}(k(\varepsilon_n)+a)\,\,  \in \,\, ]0, \varepsilon_1^\star[.
\end{equation}
 It is immediate to check that:
\begin{equation}
\label{sequence2}
k\left(\varepsilon_{(a,n)}\right)= \delta_1 \ln (\varepsilon_{n})+a=a \pmod{2\pi}.
\end{equation}

The following proposition establishes  the convergence of the map $\mathcal{T}_{\varepsilon_{(a,n)}}$  to a singular limit as $n \rightarrow +\infty$, ($\|\star\|_{\textbf{C}^r}$ represents the norm   in the $C^r$--topology for $r\geq2$):

\begin{lemma}
\label{important lemma}
The following equality holds:
$$
\lim_{n\in \NN} \|\mathcal{T}_{\varepsilon_{(a,n)}} (x,\overline{y}) -(h_a(x,\overline{y}), \textbf{0})\|_{\textbf{C}^3} =0$$
where $\textbf{0}$ is the constant null map and
\begin{equation}
\label{circle map}
h_a(x, \overline{y})= x+\alpha_1+a+\delta_1 \ln (\overline{y} +\Psi_2(x,\overline{y})) 
\end{equation}
\end{lemma}

\begin{proof}

Using \eqref{sequence2}, at $\varepsilon_1={\varepsilon_{(a,n)}}$, we get:
{ \begin{eqnarray*}
\mathcal{F}_1 (x,\overline{y})&=&x+\alpha_1+{\varepsilon_{(a,n)}}  \Psi_1(x,\overline{y})+\delta_1  \ln \varepsilon_{(a,n)} +\delta_1 \ln  \left[\left(\overline{y}+\Psi_2(x,\overline{y})\right)\right] \pmod{2\pi} \\
&=&x+\alpha_1+{\varepsilon_{(a,n)}}  \Psi_1(x,\overline{y})+a +\delta_1\ln  \left[\left(\overline{y}+\Psi_2(x,\overline{y})\right)\right]\pmod{2\pi}
  \end{eqnarray*}}and
  \begin{eqnarray*}
\mathcal{F}_2 (x,\overline{y}) &=& {\varepsilon_{(a,n)}}^{\delta-1} \left(\overline{y}+ g(x,\overline{y})\right)^\delta.
  \end{eqnarray*}
    Therefore, since $\dpt \lim_{n\in \NN} {\varepsilon_{(n,a)}}=0$ we may write:
  $$
\lim_{n\in \NN}   \mathcal{F}_1 (x,\overline{y}) = x+\alpha_1+a+ \delta_1 \ln  \left[\left(\Psi_2(x,0)\right)\right] \pmod{2\pi}
  $$
and
  $$
\lim_{n\in \NN}   \mathcal{F}^2_{\varepsilon_{(a,n)}} (x,\overline{y}) = 0,
  $$
 and we get the result (here we make use of  the condition $\delta>1$ in \textbf{(H2)}).\\
\end{proof}

\begin{figure}[h]
\begin{center}
\includegraphics[height=5cm]{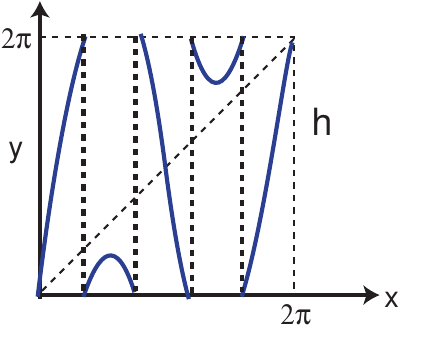}
\end{center}
\caption{\small   Effects of large $\delta_1$ on the graph of a  Misiurewicz-type map $h:\EU^1 \rightarrow \EU^1$. Illustration of the mixing property \textbf{(P7)}.  }
\label{Bruno1.2}
\end{figure}

\begin{remark}
\label{rem7.2}
The map $h_a(x) =x+\alpha_1+a  +\delta_1 \ln  \Psi_2(x,0)\equiv   \mathcal{F}^1_{\varepsilon_{(a,n)}} (x,0)$ has two nondegenerate critical points, by Hypothesis \textbf{(H4)}. See in Figure \ref{Bruno1.2} the effect of large $\delta_1$ on the graph of $h_a$.\\
\end{remark}

\subsection{The proof}
Using the terminology introduced in Subsection \ref{s: limit family}, we can conclude that:
\begin{lemma}
\label{tec.lemma}
The family $F_{(a,b)}=\mathcal{T}_{(\varepsilon_1,0)}$ for $b=\varepsilon_{(a,n)}$ satisfies \textbf{(P1)--(P7)} itemised in  Subsection \ref{rank_one}. 
\end{lemma}

\begin{proof}
Property \textbf{(P1)} follows from the same reasoning of \cite[\S 8.4]{Rodrigues2022}; \textbf{(P2)--(P3)} follow from Lemma \ref{important lemma}; \textbf{(P4)--(P5)} and \textbf{(P7)} are a consequence of Proposition \ref{Prop2.1WY} combined with Lemma \ref{important lemma};  \textbf{(P6)} follows from Remark \ref{rem7.2}. \\
\end{proof}

From Lemma \ref{tec.lemma}, we may conclude that, for $\varepsilon_1<\varepsilon_1^\star$ and $\delta_1 \gg1$,  there exists $ \mathcal{G}_1 \subset  [0, \varepsilon_1^\star]$ with positive Lebesgue
measure such that if $\varepsilon_1 \in \mathcal{G}_1$, then the map $\mathcal{T}_{(\varepsilon_1, 0)}$ admits an irreducible strange attractor $\Omega$ such that
\begin{equation}
\label{Omega}
 \Omega \subset \bigcap_{m=0}^{+\infty}  \mathcal{T}_{(\varepsilon_1, 0)}^m(\mathcal{A})
 \end{equation}
shadowing $ y=1\,  (\Leftrightarrow\overline{y}=0)$ and supporting a unique ergodic SRB measure $\nu$. The orbit of Lebesgue almost all points in ${\Omega}$ has positive Lyapunov exponent  and is asymptotically distributed according to $\nu$.

\begin{remark}
Since the map $h_{a^\star}$ has a periodic point $p_a'$ belonging to a transverse homoclinic orbit, then for $\varepsilon_1$--sufficiently small and for $a$ close to $a^\star$, any $\varepsilon_1$--perturbation of $(h_{a^\star}(x,\textbf{0}), \textbf{0})$ (in coordinates $(x, \overline{y})$) possesses a hyperbolic  periodic point (of saddle type) $p_a$ which is the analytic continuation of $(p_a',\textbf{0})$. This follows from \cite{WY}.
\end{remark}
 
\section{Proof of Theorem \ref{prop0}}
    \label{s: Th B proof} 
 
Let us define $p,q\in \NN$ as in Section \ref{s: main}. Fix $\varepsilon_1\in \mathcal{G}_1$ (from Theorem \ref{th:transv}) and $\varepsilon_2=0$.\\

Note that if     $( \alpha_2, \delta_2) \notin \mathcal{U}^{p/q}$, then the map $\mathcal{F}_3$ is conjugated to an irrational rotation on the circle, and $\mathcal{T}_{(\varepsilon_1,0)}$ and $\mathcal{F}_3$ are uncoupled.\\
 
  \begin{enumerate}
  \item Proving that $\mathcal{F}_{(\varepsilon_1, 0)}= \mathcal{T}_{(\varepsilon_1, 0)}\times \mathcal{F}_3$ has a dense orbit in $\Omega \times \EU^1$   is equivalent to   show that   for any non-empty open sets $U', V' \in \Omega \times \EU^1$, there exists $k\in \NN$ such that $\mathcal{F}_{(\varepsilon_1, 0)}^k(U')\cap V'\neq \emptyset$.
Take the following non-empty open sets of $\Omega\times \EU^1$ (see Figure \ref{thB proof}):
$$
U'=U\times [r-\tau_1, r+\tau_1], \quad \text{and}\quad V'=V\times [s-\tau_2, s+\tau_2],
$$
where:
\begin{itemize}
\item $\tau_1, \tau_2\gtrsim 0$, 
\item $r, s \in [0,2\pi]$ and 
\item $U,V$ are non-empty open sets of $\Omega \subset \EU^1\times [1, 1+b]$, endowed with the induced topology.\\
\end{itemize}

For a fixed $\tau_2>0$, we can find a sequence $(n_j)_{j\in \NN}$ such that $\mathcal{F}_3^{n_j}(r)\in  [s-\tau_2, s+\tau_2]$, where 
\begin{equation}
\label{seq1}
0<n_1<N \quad \text{and} \quad n_{j+1}-n_j<N
\end{equation}
for all $j\in \NN$, and  $N$ is independent on  $r$ and $s$.\\

 \begin{figure}[h]
\begin{center}
\includegraphics[height=5cm]{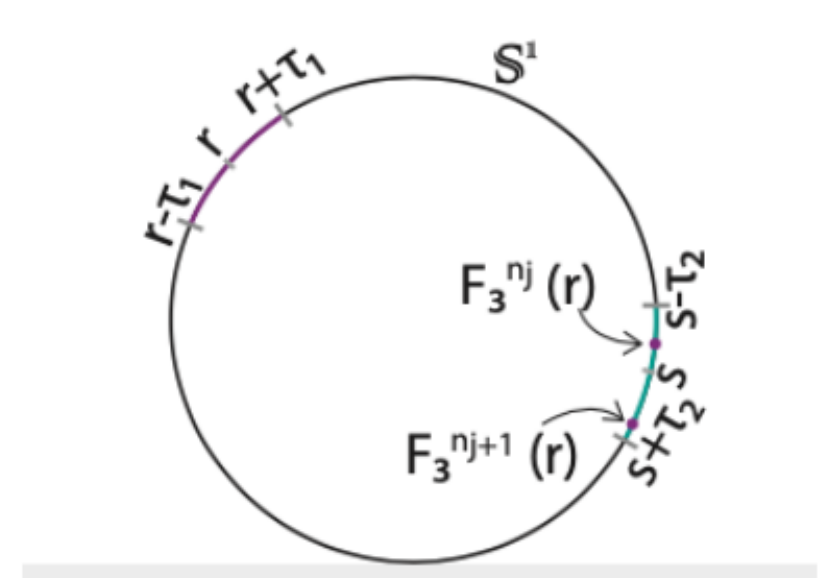}
\end{center}
\caption{\small   Illustration of the sets  $ [r-\tau_1, r+\tau_1]$ and  $[s-\tau_2, s+\tau_2]$ in $\EU^1$ and possible location of the sequence $\mathcal{F}_3^{n_j}(r)$.}
\label{thB proof}
\end{figure}

By Theorem \ref{th:transv}, we know that there exists a hyperbolic periodic orbit $p_a$ such that $ W^u(p_a)$ is dense in  $\Omega$ (since  $\overline{W^u(p_a)}=\Omega$ and $\varepsilon_1\in \mathcal{G}_1$). Then,  there exists a point $q' \in W^u(p_a)\cap V$. Consider $u=\mathcal{T}_{(\varepsilon_1, 0)}^{-l}(q')$, $l\gg N$, such that $u$ and its first $N$ iterates are arbitrarily close to $p_a$ (with respect to the metric ${\rm dist}_2$).
Then,  there are $N+1$ non-empty open sets $Z_0$, $Z_1$,..., $Z_N$ centered at $u,\mathcal{T}_{(\varepsilon_1, 0)}(u), ..., \mathcal{T}_{(\varepsilon_1, 0)}^N(u)$ such that their images under $\mathcal{T}_{(\varepsilon_1, 0)}^l, \mathcal{T}_{(\varepsilon_1, 0)}^{l-1}, ..., \mathcal{T}_{(\varepsilon_1, 0)}^{l-N}$ are contained in $V$.

Shrinking (if necessary) the neighbourhood of $p_a$, there exists a point $x\in U\cap W^u(p_a)$ belonging to a dense orbit and having positive Lyapunov exponent, such that $\mathcal{T}_{(\varepsilon_1, 0)}^m(x)\in Z_j$, for some $m\in \NN$ and all $j\in \{0, ..., N\}$.  In particular,  
$$
\mathcal{T}_{(\varepsilon_1, 0)}^{l+m-j}(U)\cap V\neq\emptyset
$$
for all $j=0,..., N$. 
Observing that 
$$
l+m-N\leq l+m-j\leq l+m,
$$
then there exists $j\in \{0, ..., N\}$ such that $n_i<l+m-j<n_{i+1}$, for some $n_i$ defined in \eqref{seq1}. This means that some iterates of $\mathcal{F}_3^{l+m-j}(r)$, $j\in \{0, ..., N\}$, lie inside the interval $(s-\tau_2, s+\tau_2)$, which implies  that there exists $j\leq N$ such that $l+m-j=n_i$ for which:
$$
\mathcal{F}_{(\varepsilon_1, 0)}^{l+m-j}(U')\cap V'\neq \emptyset.
$$

\item The existence of   $k>0, \lambda>1$, $z=(x_0,y_0, t_0)\in \mathcal{A}'$ and a vector $v=(v_x, v_y, v_t)\in T_zM$ such that:\\
  $$\left\|D\mathcal{F}_{(\varepsilon_1,0)}^n(z)v\right\|\geq k\lambda ^n,\quad \forall n\in \NN_0,$$
  follows from  Theorem \ref{th:transv} combined with the inequality  (\footnote{Observe that the domain of the norm of the left hand side of \eqref{ineq1} is $M$ and, on the right hand side, is  $\EU^1\times [1, 1+b]$.}):
\begin{equation}
\label{ineq1}
  \left\|D\mathcal{F}_{(\varepsilon_1,0)}^n(z)v\right\| \geq \  \left\|D\mathcal{T}_{(\varepsilon_1,0)}^n(x_0, y_0)(v_x, v_y)\right\|
  \end{equation}
  where $n\in \NN$ and  $v, z$ are as above.
\end{enumerate}

  \section{Proof of Theorem \ref{Th2.1}}
  \label{s: Prova Th C}
   The aim of this section is the proof of Theorem \ref{Th2.1}. We start by reviewing the Tangerman-Szewc Theorem  which will be useful in the sequel. 
  \begin{figure}[h]
\begin{center}
\includegraphics[height=3.5cm]{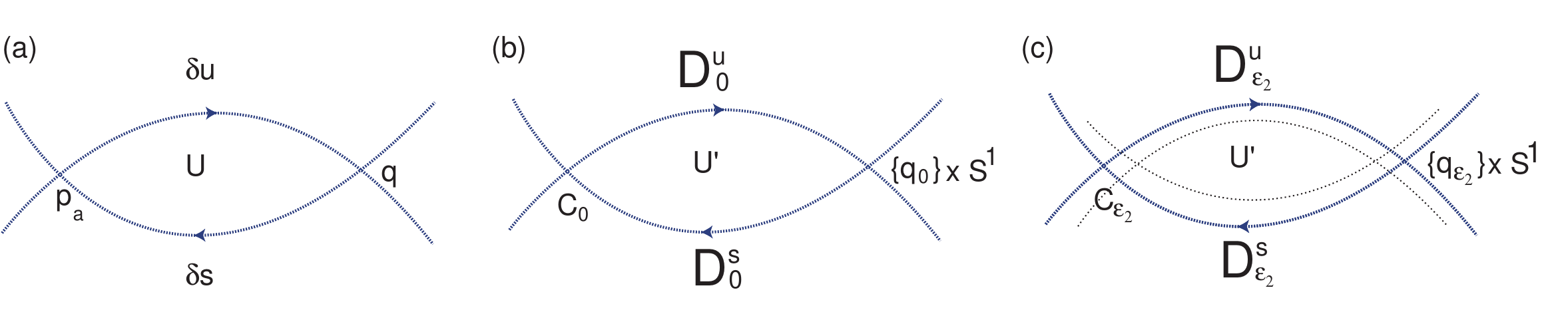}
\end{center}
\caption{\small  (a):  Homoclinic orbit of  $p_a\in \EU^1\times [1, 1+b]$ associated to $\mathcal{T}_{(\varepsilon_1, 0)}$. (b)   Homoclinic orbit of  $\mathcal{C}_0 $ associated to $\mathcal{F}_{(\varepsilon_1, 0)}$.  (c)   Homoclinic orbit of  $\mathcal{C}_{\varepsilon_2}   $ associated to $\mathcal{F}_{(\varepsilon_1, \varepsilon_2)}$. Dashed lines in (c) corresponds to the Homoclinic orbit of  $ \mathcal{C}_0 $.  }
\label{Bruno1.3}
\end{figure}
  \subsection{Tangerman-Szewc Theorem revisited}
\label{ss:tangerman}
Let $K\subset \EU^1\times \RR$ be a compact set. Let $$T : K \to  \EU^1\times \RR$$ be a diffeomorphism having a   dissipative saddle fixed point $p $ such that:
 
\begin{itemize}
\item the stable and unstable manifolds $W^s(p)$ and $W^u(p)$ intersect transversally at the homoclinic point $q \in W^s(p) \cap W^u(p)$ and
\item $W^u(p)$ is bounded as a subset of $K\subset  \EU^1\times \RR$.
 \end{itemize}  
 
The F. Tangerman and B. Szewc Theorem \cite[Appendix 3]{PT} states that there is a nonempty open set $U\subset K$ such that for each $p\in U$, $\omega(p)\subset \overline{W^u(p)}$.  As suggested by Figure \ref{Bruno1.3}(a), this open set, say $U$, may be taken as the region bounded by the two arcs, say $\delta s \subset W^s(p)$ and  $\delta u \subset W^u(p)$. \\

\subsection{Proof of  Theorem \ref{Th2.1}(1) for $\varepsilon_2=0$ (uncoupled system)}
\label{provaC(1)}
\,\\ 

Fix $\varepsilon_1 \in [0, \varepsilon_1^\star]$ and consider the circle $\mathcal{C}_0 =\{p_a\}\times \EU^1 $, invariant under the map $\mathcal{F}_{(\varepsilon_1, 0)}$, where $p_a$ is the hyperbolic periodic orbit given by Theorem \ref{th:transv}.
 The manifolds $W^u(\mathcal{C}_0)$ and $W^s(\mathcal{C}_0)$ are given by $W^u(p_a)\times \EU^1$ and $W^s(p_a)\times \EU^1$, respectively. As shown in Figure \ref{Bruno1.3}(b), they intersect transversally at a circle $\{q_0\} \times \EU^1$, consisting of points homoclinic to $\mathcal{C}_0 $.

As before, consider the two arcs $\delta s$ and $\delta u$ with extremes $p_a\in \EU^1\times [1,1+b]$ and $q_0\in \EU^1\times [1,1+b]$, respectively. They limit an open set $U\subset \EU^1\times [1,1+b]$.
 Define $D^s_0$ and $D^u_0$ to be the portions of stable, and unstable manifold of $\mathcal{C}_0$, respectively, given by
\begin{equation}
\label{Ds}
D^s_0 =\delta s \times \EU^1\quad \text{and}\quad  D^u_0= \delta u \times \EU^1.
\end{equation}

  Both sets $D^s_0$ and $D^u_0$ are compact, and their union forms the boundary of the open region $U' = U \times \EU^1\subset M$, which is a  topological  solid torus (since $p_a\neq (0,0)$).\\

  The following technical result is valid for $U'\subset M$ and  says that  the forward evolution of every point $(x,y, t)\in U'$ approaches the boundary of $\mathcal{F}_{(\varepsilon_1,0)}^n(U')$,  denoted by  $\partial \mathcal{F}^n_{(\varepsilon_1,0)}$, as $n \to +\infty$. 
\begin{lemma} The following assertions hold for  $(x,y, t)\in U'\subset M$. \\
\begin{enumerate}
\item $\dpt \lim_{n\in \NN}  \ell_3(\mathcal{F}_{(\varepsilon_1,0)}^n(U'))=0$ and \\
\item $\dpt \lim_{n\in \NN}{\rm dist_3}(\mathcal{F}_{(\varepsilon_1,0)}^n(x,y,t), \partial \mathcal{F}_{(\varepsilon_1,0)}^n(U'))=0$.
\end{enumerate}

\end{lemma}
\begin{proof}
\begin{enumerate}
\item This item follows from the chain of equalities:
\begin{eqnarray}
\nonumber \ell_3(\mathcal{F}_{(\varepsilon_1,0)}^n(U')) &=& 2\pi \int_{\mathcal{T}_{(\varepsilon_1,0)}^n(U)}\, {\rm dx dy} \\
\nonumber && \\
 \nonumber& {\leq }&   2\pi \int_{U} |\det D\mathcal{T}_{(\varepsilon_1 0)}^n(x,y)|{\rm dx dy}\\
 \nonumber && \\
 \nonumber &\overset{\text{Remark}\,  \eqref{cont1}}\leq &2\pi K^n \ell_2(U) , \quad \text{for some} \quad K\in (0,1).
 \end{eqnarray}
Since $U$ and $U'$ are bounded and
 $$
0\leq   \lim_{n\in \NN}  \ell_3(\mathcal{F}_{(\varepsilon_1,0)}^n(U')) \leq  \lim_{n\in \NN} 2\pi K^n \ell_2(U)=0,
 $$
 the result follows. \\
 \item Suppose, by contradiction, that item (2)  does not hold. Then, 
there exists $c>0$ such that for all $n\in \NN$, there exists $n_0\in \NN$ such that 
$$
n\geq n_0\quad \text{and}\quad {\rm dist_3}\left(\mathcal{F}_{(\varepsilon_1,0)}^n(x,y,t), \partial \mathcal{F}_{(\varepsilon_1,0)}^n(U')\right)>c.
$$
 In particular, for $n\in \NN$, the open ball centered at $\mathcal{F}_{(\varepsilon_1,0)}^n(x,y,t)$ and radius $c>0$ contained in $\mathcal{F}_{(\varepsilon_1,0)}^n(U'),$ for all $n$. This would contradict the first item. Then the result follows.

 \end{enumerate}
\end{proof}

Observing \eqref{Ds}, the boundary of $\mathcal{F}_{(\varepsilon_1,0)}^n(U')$   consists of two portions of stable and unstable manifold of $\mathcal{C}_0$:
$$
\partial \mathcal{F}_{(\varepsilon_1,0)}^n(U')= \mathcal{F}_{(\varepsilon_1,0)}^n(D^s_0)\cup \mathcal{F}_{(\varepsilon_1,0)}^n(D^u_0).\\
$$
\\

Due to  $\dpt \lim_{n\in \NN}{\rm diam}(\mathcal{F}_{(\varepsilon_1,0)}^n(D_0^s)) =0$ (all points in $D^s_0$ go to $\mathcal{C}_0$, by definition of stable manifold of a hyperbolic set) and since $W^u(\mathcal{C}_0)$ is bounded, then all points starting at $U'\backslash W^s(\mathcal{C}_0)$ are bounded and approach $W^u(\mathcal{C}_0)$. Formally, we may write:
$$
\forall (x,y,t)\in U', \quad \lim_{n\in \NN}{\rm dist_3} \left(\mathcal{F}_{(\varepsilon_1,0)}^n(x,y, t), \mathcal{F}_{(\varepsilon_1,0)}^n(D^u_0)\right)=0
$$
which implies: 
$$
\forall (x,y,t)\in U', \quad \omega(x,y,t) \subset \overline{W^u(\mathcal{C}_0)}\\
$$
and    Theorem \ref{Th2.1}(1) is proved for $\varepsilon_2=0$.
\bigbreak
\subsection{Proof of  Theorem \ref{Th2.1}(1) for $\varepsilon_2\gtrsim0$ (coupled system)}

Let $\varepsilon_2<\varepsilon_2^\star$ where $\varepsilon_2^\star$ is the value given by the persistence of normally hyperbolic invariant manifolds (in the $C^3$--topology). Then $\mathcal{F}_{(\varepsilon_1,\varepsilon_2)}$   has a   hyperbolic invariant circle $\mathcal{C}_{\varepsilon_2}$ of saddle type and the manifolds $W^u(\mathcal{C}_{\varepsilon_2})$ and $W^s(\mathcal{C}_{\varepsilon_2})$  are $C^3$--close to $W^u(\mathcal{C}_0)$ and $W^s(\mathcal{C}_0)$. We assume that the latter manifolds intersect transversally at $\{q_0 \}\times \EU^1$.\\

\begin{lemma}
\label{transv}
For $\varepsilon_2<\varepsilon_2^\star$, the invariant manifolds $W^u(\mathcal{C}_{\varepsilon_2})$ and $W^s(\mathcal{C}_{\varepsilon_2})$  intersect transversally. 
\end{lemma}
\begin{proof}
We are going to use the classical theory of hyperbolic manifolds (see for instance \cite{PM}).
First, consider the arcs $\delta u \subset W^u(p_a)$ and  $\delta s \subset W^s(p_a)$ as in Figure \ref{Bruno1.3}(a). Define
$$
A^u_0=\delta u\times \EU^1, \quad A^s_0= \delta s\times \EU^1.
$$
The circle $\{q_0\}\times \EU^1$ is the intersection of the manifolds $A^u_0$ and $A^s_0$, bounded away from their frontiers.
Now, consider the inclusions $$i_0 : A^u_0 \to M\quad \text{and}\quad j_0 : A^s_0 \to M.$$
     By the closeness of $W^u(\mathcal{C}_0)$ and $W^u(\mathcal{C}_{\varepsilon_2})$, there exists a
$C^3$--diffeomorphism 
$$h: A^u_0 \to A^u_{\varepsilon_2}\subset W^u(\mathcal{C}_{\varepsilon_2})$$ 
 such that the map $i_0$ is $C^3$--close to $i_{\varepsilon_2} \circ h$, where $ i_{\varepsilon_2} : A^u_{\varepsilon_2} \to M$ is the inclusion  map.  
  Analogously, there exists a $C^3$--diffeomorphism 
$$k: A^s_0 \to A^s_{\varepsilon_2}\subset W^s(\mathcal{C}_{\varepsilon_2})$$ 
 such that the map $j_0$ is $C^3$--close to $j_{\varepsilon_2} \circ k$, where $ j_{\varepsilon_2} : A^s_{\varepsilon_2} \to M$ is again the inclusion map.

 The map 
  $$
  i_0\times j_0: A^u_0\times A^s_0 \to M^2$$ is transversal to the diagonal $\Delta=\{(x,x),x\in M\}$ since the manifolds  $W^u(\mathcal{C}_{0})$ and $W^s(\mathcal{C}_{0})$  intersect transversally.   For $\varepsilon_2\gtrsim 0$ small, the map $$(i_{\varepsilon_2} \circ h)\times (j_{\varepsilon _2}\circ k): A^u_0 \times A^s_0 \to  M^2$$ is $C^3$--close to $i_0\times j_0$.
     Since $\Delta$ is closed and $A^u_0 \times A^s_0$ is compact,  then there exists $\varepsilon_2^{\star\star}<\varepsilon_2^\star$ such that $(i_{\varepsilon_2} \circ h)\times (j_{\varepsilon_2} \circ k) \pitchfork \Delta$ for $\varepsilon_2<\varepsilon_2^{\star\star}$, and  the submanifolds
$$
(i_0\times j_0)^{-1}(\Delta) \quad \text{and} \quad 
(i_{\varepsilon_2} \circ h)  \times (j_{\varepsilon_2} \circ k) ^{-1}(\Delta) 
$$
are diffeomorphic. Since $[(i_{\varepsilon_2} \circ h)  \times (j_{\varepsilon_2} \circ k)] ^{-1}(\Delta) $ is diffeomorphic to $A^u_{\varepsilon_2} \cap A^s_{\varepsilon_2} $ for all $\varepsilon_2<\varepsilon_2^{\star\star}$, and 
$$
(i_0\times j_0)^{-1}(\Delta) = A^u_0\cap A^s_0= \{q_0\}\times \EU^1,
$$
then   the intersection $ A^u_{\varepsilon_2} \cap  A^s_{\varepsilon_2} $ is diffeomorphic to $  \{q_{\varepsilon_2}\}\times \EU^1$, as suggested by Figure \ref{Bruno1.3}(c).
  \end{proof}

Now, we prove the result by using a similar argument to that of  Subsection \ref{provaC(1)}. Indeed, define $D_{\varepsilon_2}^u$ as the compact part of $W^u(\mathcal{C}_{\varepsilon_2})$ bounded by the invariant circle $\mathcal{C}_{\varepsilon_2}$ and the circle of homoclinic points $\{q_{\varepsilon_2}\}\times \EU^1$. Define $D_{\varepsilon_2}^s$ in an analogous way.

The manifolds $D_{\varepsilon_2}^u\subset W^u(\mathcal{C}_{\varepsilon_2})$ and  $D_{\varepsilon_2}^s\subset W^s(\mathcal{C}_{\varepsilon_2})$ form the boundary of an open region $U'\subset M$ homeomorphic to a torus (note that the set $U'$ may be not the same as that of Subsection \ref{provaC(1)}).
    By the $C^3$--closeness of the perturbed manifolds $W^u(\mathcal{C}_{\varepsilon_2})$ and $W^s(\mathcal{C}_{\varepsilon_2})$ to the unperturbed $W^u(\mathcal{C}_0)$ and $W^s(\mathcal{C}_0)$, both $U'$ and $W^u(\mathcal{C}_{\varepsilon_2})$  
  are bounded.\\

The map $\mathcal{F}_{(\varepsilon_1,\varepsilon_2)}$   is dissipative (cf. Remark \ref{cont1}). Indeed, by taking $\varepsilon_2<\varepsilon_2^{\star\star}$ small enough, we ensure that 
   $$
   |\det(D\mathcal{F}_{(\varepsilon_1,\varepsilon_2)}(x,y,t))| <  1
   $$ 
   for all   $(x,y,t) \in U'$.  Like in the first part of the proof performed in  Subsection \ref{provaC(1)}, one has
$$
\forall (x,y,t)\in U', \quad \omega(x,y,t) \subset \overline{W^u(\mathcal{C}_{\varepsilon_2})}
$$
and we get the result  (by taking $\varepsilon_2^{\star\star}=\varepsilon_2^\star$).

 \subsection{Proof of  Theorem \ref{Th2.1}(2)}
The existence of hyperbolic horseshoes follows from Lemma \ref{transv} and \cite{PT}. The transverse intersection of invariant manifolds of a hyperbolic saddle periodic orbit  is a sufficient condition for the existence of horseshoes whose non-wandering set is dynamically defined by $\overline{W^u(\mathcal{C}_{ \varepsilon_2})\pitchfork {W^s(\mathcal{C}_{ \varepsilon_2}})}$.

\subsection{Proof of Theorem \ref{Th2.1}(3)}
  This proof runs along the same lines to those Proposition 2.7 of \cite{BST98}. 
  
Let $p_a = (x_0, y_0)\in \EU^1\times [1,1+b]$ be a saddle
fixed point of the diffeomorphism $\mathcal{T}_{(\varepsilon_1,0)}$. 
The invariant circle $\mathcal{C}_{0} = \{p_a\} \times \EU^1$ of $\mathcal{F}_{(\varepsilon_1,0)}$
can be  seen as a graph over $\EU^1$:
$$
\mathcal{C}_0=\{(x_0, y_0, t),t\in \EU^1\}.
$$
Fix   $\varepsilon_2 < \varepsilon_2^\star$(\footnote{The value $\varepsilon_2^\star$ is given by the persistence of normally hyperbolic invariant manifolds (in the $C^3$--topology).}), where $\varepsilon_2^\star$ is taken from   Theorem \ref{Th2.1}(1). By the $C^3$--closeness of $\mathcal{C}_0$ and $\mathcal{C}_{\varepsilon_2}$, the invariant circle $\mathcal{C}_{\varepsilon_2}$ of  $\mathcal{F}_{(\varepsilon_1,\varepsilon_2)}$ can be written as a $C^3$--graph over $\EU^1$ as follows:
$$
\mathcal{C}_{\varepsilon_2}=\{(x_{\varepsilon_2}(t), y_{\varepsilon_2}(t), t),t\in \EU^1\}
$$
where $x_{\varepsilon_2}, y_{\varepsilon_2}:\EU^1 \to \RR$ are of the type
\begin{eqnarray*}
x_{\varepsilon_2}(t)&=& x_0 +\mathcal{O}(\varepsilon_2)\\
y_{\varepsilon_2}(t)&=& y_0 +\mathcal{O}(\varepsilon_2),
  \end{eqnarray*}
where $\mathcal{O}$ stands for the Landau notation. Therefore, the restriction of
 $\mathcal{F}_{(\varepsilon_1,\varepsilon_2)}$  to $\mathcal{C}_{\varepsilon_2}$  can be seen as a map on $\EU^1$.
Doing the Taylor expansion in $\varepsilon_2$, there exists a smooth map $c: \RR^2\to \RR$ depending on $\alpha_2, \varepsilon_2\in \RR$, such that $Y_{\varepsilon_2}(t)$  is conjugate to the normal form:
$$
NF(Y_{\varepsilon_2}(t)))= t+\alpha_2+c(\alpha_2, \varepsilon_2)+\mathcal{O}(\varepsilon_2 ^{r+1}),
$$ 
for $r>1$.  The conjugacy is explicitly given in Sections 2.2 and 2.3 of \cite{BST98}.
Fixing  $\tau = 3$ and $\gamma > 0$ in the diophantine condition of the proof of Proposition 2.7 of \cite{BST98}, there exists  a set $\Gamma_{\varepsilon_2}$ with 
positive Lebesgue measure such that: if $\varepsilon_2\in \Gamma_{\varepsilon_2}$, then  the family $\mathcal{F}_{(\varepsilon_1,\varepsilon_2)}$, restricted to $\mathcal{C}_{\varepsilon_2}$,
is conjugate to a family of irrational rotations. The conjugacy is, at least, of class $C^5$ (cf. \cite{BST98}).  This proves   the Theorem.

  \section{Proof of Theorem \ref{Th2.5}}
  
 \label{s: Prova Th D}
 
We divide the proof into two cases, according to the rotation number of the Arnold tongue under consideration. We remind the reader that the value of  $\delta_2$  satisfies Hypothesis \textbf{(H5)}  -- see Remark \ref{injective}.
  
\subsection{The case on the rotation number $1$}
 \label{rot_num1}
 Fix $\delta_1\gg1 $ and $\varepsilon_1 \in \mathcal{G}_1$ verifying the hypotheses of Theorem \ref{th:transv}.  This subsection considers the case $(\alpha_2, \delta_2) \in {\rm int} (\mathcal{U}^1)$, the interior of a resonance tongue of rotation number one -- cf. Figure \ref{Bruno1.1}.  Then the map $\mathcal{F}_3$ on $\EU^1$ has two hyperbolic fixed points $t_1^s$ (attracting) and $t_1^r$ (repelling), where the $t$--coordinate of both points depends on the choice of $(\alpha_2, \delta_2) \in {\rm int}(\mathcal{U}^1)$. 
For all $t \in \EU^1$ with $t\neq t^r_1$, the orbit of $t$ under $\mathcal{F}_3$ converges to $t^s_1$. This means that the manifold
  $$
\Omega_1=\{(x,y, t)\in \EU^1\times [1,1+b] \times \EU^1: t= t_1^s\}
  $$
  is $\mathcal{F}_3$--invariant and globally attracting in $\mathcal{A}'$. Define the map $G_{1}:\Omega_1 \to \Omega_1$ the restriction of $\mathcal{F}_{(\varepsilon_1, \varepsilon_2)}$
to  $\Omega_1$ whose expression is given by:
\begin{equation*}
\left(\begin{array}{c}x\\ \\ y\\ \\t_1^s\end{array}\right)
\mapsto
 \left(\begin{array}{c}
x+ \alpha_1 + \varepsilon_1 \Psi_1(x,y) + \delta_1\ln ( (y-1)+ \varepsilon_1[ \Psi_2 (x,y)+\varepsilon_2\Psi_4 (x,y,t_1^s)]) \pmod{2\pi} 
\\
\\
1+ ((y-1)+ \varepsilon_1 g(x,y))^\delta
\\
\\
t_1^s 
\end{array}\right).
\end{equation*}

The map $G_{1}$ is well defined.
    By construction, one knows that
        \begin{eqnarray}
    \label{basin1}
\nonumber G_{1}(\mathcal{A}') &\subset &\mathcal{A}', \\
\nonumber G_{1}(\mathcal{A} \times \{t_1^s\})&\subset &{\rm int} (\mathcal{A} \times \{t_1^s\}), \\
G_{1}(\mathcal{A}  \times (\EU^1\backslash \{t_1^r\}))&\subset & {\rm int} (\mathcal{A}  \times (\EU^1\backslash \{t_1^r\})).
  \end{eqnarray}

    \bigbreak
    Since the set $\Omega_1$ is diffeomorphic to $\RR^2$, we may consider $G_1$ as a map of $\EU^1\times [1,b+1]$. Then, in the appropriate domain,  $G_1$ is
a $  \mathcal{O}(\varepsilon_2)$--perturbation of the map $\mathcal{T}_{(\varepsilon_1, 0)}$ and its singular limit is of the form
$(\overline{x}, \overline{y})=(\tilde{h}_a(x), 0)$, where $\|\tilde{h}_a-{h}_a\|_{\textbf{C}^3}=  \mathcal{O}(\varepsilon_2)$.

 Using Proposition \ref{Prop2.1WY}, $\mathcal{G}_{1}$  satisfies \textbf{(P1)--(P7)} of  \cite{WY}  and then  we may conclude that for 
$\delta_1^\star> 0$ and $ 0<\varepsilon_1^\star\ll 1$  there exists $\tilde{\mathcal{G}}_1\subset [0, \varepsilon_1^\star]$ with positive Lebesgue measure for every $\varepsilon_1\in \tilde{\mathcal{G}}_1$ the map $G_1$  exhibits an irreducible strange attractor $\tilde{\Omega}$ that supports a unique ergodic SRB measure $\nu$.  Furthermore, we have $\tilde{\Omega}\subset \overline{W^u{(p_a)}}$, where $p_a$ is a periodic orbit of $G_1$ and  shadows the entire circle defined by $\EU^1\times \{1\} $ (in coordinates $(x,y)$), as a consequence of   Theorem \ref{th:transv}.\\

Now we need to ``transfer'' the dynamical information from $G_1$ to the map $\mathcal{F}_{(\varepsilon_1, \varepsilon_2)}$. \\

The point $p^\star = (p_a,t_1^s)=(x_0, y_0, t_1^s)$ is a saddle periodic point of the map $\mathcal{F}_{(\varepsilon_1, \varepsilon_2)}$, and $$W^u(p^\star)=W^u(p_a)\times \{t_1^s\}.$$ Therefore
$\overline{W^u(p^\star)}= \tilde{\Omega}\times \{t_1^s\}$ is a rank-one strange attractor of  $\mathcal{F}_{(\varepsilon_1, \varepsilon_2)}$ because:
\begin{itemize}
\item it has a dense orbit  $\tilde{z}=(\tilde{x}, \tilde{y}, t_1^s)$ with  a positive Lyapunov exponent. This comes from    Theorem \ref{th:transv} and
  $$
  \left\|D\mathcal{F}_{(\varepsilon_1,0)}^n(\tilde{z})(v_x, v_y, v_t)\right\| \geq \  \left\|D{G}_1^n(\tilde{x}, \tilde{y})(v_x, v_y)\right\|>k\lambda^n,
  $$
  where $k>0$ and $\lambda>1$;
 
\item $\overline{W^u(p^\star)}$ has nonempty interior in $\mathcal{A}\times \EU^1$  because of (\ref{basin1}).
\end{itemize} 

 By construction the set $\overline{W^u(p^\star)}$ shadows the circle $\EU^1\times \{1\} \times \{t_1^s\}$ in the original coordinates $(x,y,t)\in M$.   This proves Theorem \ref{Th2.5} for $(\alpha_2, \delta_2) \in  {\rm int} (\mathcal{U}^1)$ .

  \subsection{The case on the rotational number $p/q$}
    Fix $\delta_1\gg1 $ and $\varepsilon_1 \in \mathcal{G}_1$ verifying the hypotheses of Theorem \ref{th:transv}.  Here, we consider the case $(\alpha_2, \delta_2) \in {\rm int} (\mathcal{U}^{p/q})$, the interior of the tongue of period $p/q$.
   Then $\mathcal{F}_3$ has a even number of repelling and attracting orbits whose orbits can be written as
 $$
 {\rm Orb}_{\mathcal{F}_3}(t_1^s)= \{t_1^s, t_2^s,..., t_q^s\} \quad \text{and} \quad  {\rm Orb}_{\mathcal{F}_3}(t_1^u)= \{t_1^u, t_2^u,..., t_q^u\}
 $$
 
  As before, for $j\in \{1, ..., q\}$, define the union of $q$ smooth manifolds
     $$
\Omega_j=\{(x,y, t)\in \EU^1\times[1,1+b]\times \EU^1: t= t_j^s\}
  $$
  and set the maps $G_j$ as the restriction of $\mathcal{F}_{(\varepsilon_1, \varepsilon_2)}$ to $\Omega_j$ as follows:
  $$
  G_j: \Omega_j \to \Omega_{j+1}, \quad \text{and} \quad  G_q: \Omega_q \to \Omega_{1},
  $$
  where $G_j$, $j\in \{1, ..., q-1\}$, is defined as:
  \begin{equation*}
\left(\begin{array}{c}x\\ \\ y\\ \\t_j^s\end{array}\right)
\mapsto
 \left(\begin{array}{c}
x+ \alpha_1 + \varepsilon_1 \Psi_1(x,y) + \delta_1\ln ( (y-1)+ \varepsilon_1[ \Psi_2 (x,y)+\varepsilon_2\Psi_4 (x,y,t_j^s)]) \pmod{2\pi} 
\\
\\
1+ ((y-1)+ \varepsilon_1 g(x,y))^\delta
\\
\\
t_{j+1}^s  
\end{array}\right)
\end{equation*}
and $G_q$ is defined as:
  \begin{equation*}
\left(\begin{array}{c}x\\ \\ y\\ \\t_q^s\end{array}\right)
\mapsto
 \left(\begin{array}{c}
x+ \alpha_1 + \varepsilon_1 \Psi_1(x,y) + \delta_1\ln ( (y-1)+ \varepsilon_1[ \Psi_2 (x,y)+\varepsilon_2\Psi_4 (x,y,t_q^s)]) \pmod{2\pi} 
\\
\\
1+ ((y-1)+ \varepsilon_1 g(x,y))^\delta
\\
\\
t_1^s 
\end{array}\right).
\end{equation*}

  It is easy to check that the manifold $\Omega_1$ is invariant and attracting under $\mathcal{F}_{(\varepsilon_1, \varepsilon_2)}^q$ and for all $(x,y,t)\notin \{(x,y,t): t=t_j^r, j\in \{1,..., q\}\}$, its asymptotic dynamics is given by the map $G_q \circ G_{q-1}\circ ...\circ G_1$.
  Since $\Omega_1$ is diffeomorphic to $\RR^2$, we may consider $G_q \circ G_{q-1}\circ ...\circ G_1$ as a map of $\EU^1\times [1,1+b]$. 
  Similarly to what we did in Subsection \ref{rot_num1}, we have:\\
         \begin{eqnarray}
    \label{basin2}
 \nonumber G_q \circ G_{q-1}\circ ...\circ G_1 (\mathcal{A} \times  \{t_1^s\})&\subset &{\rm int} (\mathcal{A} \times  \{t_1^s\}, \\
G_q \circ G_{q-1}\circ ...\circ G_1(\mathcal{A}  \times\EU^1\backslash   \{t_1^r\})&\subset & {\rm int} (\mathcal{A}  \times \EU^1\backslash \{t_1^r\}).
  \end{eqnarray}

   The map $G_q \circ G_{q-1}\circ ...\circ G_1$ may be seen as a $C^3$ $\varepsilon_2$--perturbation of $\mathcal{T}_{(\varepsilon_1, 0)}$ since each $G_j$ has this property.
   Using the same line of argument used in Subsection \ref{rot_num1}, there exists $\tilde{\mathcal{G}}_1$ with positive Lebesgue measure such that the map $G_q \circ G_{q-1}\circ ...\circ G_1$ has a strange attractor which can be seen as $\overline{W^u(p_a)}$, where $p_a$ is the hyperbolic periodic orbit of $G_q \circ G_{q-1}\circ ...\circ G_1$. Observe also that $p_a$ is the hyperbolic continuation of  $p_a'$ is the saddle that generates the  strange attractor for the associated Misiurewicz map.\\
    
   Again, as before,  we need to ``transfer'' the dynamical information of $G_q \circ G_{q-1}\circ ...\circ G_1$ to the map $\mathcal{F}_{(\varepsilon_1, \varepsilon_2)}$. \\

The point $p^\star = \{p_a\} \times  {\rm Orb}_{\mathcal{F}_3}(t_1^s)$ is a saddle periodic point of the map $\mathcal{F}_{(\varepsilon_1, \varepsilon_2)}$, and $$W^u(p^\star)=W^u(p_a)\times   {\rm Orb}_{\mathcal{F}_3}(t_1^s). $$  
The set $\overline{W^u(p^\star)}$ has nonempty interior in $\mathcal{A}\times \EU^1$  because of (\ref{basin2}). The next result finishes the proof of Theorem \ref{Th2.5}.
 \begin{lemma} The map  $\mathcal{F}_{(\varepsilon_1, \varepsilon_2)}$ has 
   a dense orbit in  $\tilde{\Omega}= \overline{W^u(p^\star)}= \overline{W^u(p_a)}\times  {\rm Orb}_{\mathcal{F}_3}(t_1^s) $ with 
  a positive Lyapunov exponent.
 
\end{lemma}

\begin{proof}
 
 
Let   $(x_0, y_0, t_0)\in \tilde{\Omega}$. Then there exists $j\in \{0, 1, ..., q-1\}$ such that $\mathcal{F}_3^j(t_1^s)= t_0$.
 Since the map $G_q \circ G_{q-1}\circ ...\circ G_1$ has a strange attractor, then it has a dense orbit, say $(\tilde{x}, \tilde{y})\in \mathcal{A}$. In particular, for any $\varepsilon>0$, there exits $m\in \NN$ such that:
 $$
{\rm dist}_2( [G_q \circ G_{q-1}\circ ...\circ G_1]^m (\tilde{x}, \tilde{y}), (x_0, y_0))<\varepsilon.
 $$
 which implies that 
  $$
{\rm dist}_3\left(\mathcal{F}^{qm+j}_{(\varepsilon_1, \varepsilon_2)}(\tilde{x}, \tilde{y}, t_1^s), (x_0, y_0, t_0) \right) <\varepsilon.
$$

  If $n_1+n_2=qm+j$ and $z\in \mathcal{A}'$, by the Chain Rule, we have: 
 $$
 D \mathcal{F}_{(\varepsilon_1, \varepsilon_2)}^m(z)= D \mathcal{F}_{(\varepsilon_1, \varepsilon_2)}^{n_1}([G_q \circ G_{q-1}\circ ...\circ G_1]^{n_2}(z))\times D[G_q \circ G_{q-1}\circ ...\circ G_1]^{n_2}(z).
 $$
 
 Now, take $z=(\tilde{x}, \tilde{y}, t_1^s)\in \Omega_1$ a point having a dense orbit in $\mathcal{A}'$ and $v=(v_x, v_y, 0) \in T_z \mathcal{A}'$ such that 
 $$
 \|D[G_q \circ G_{q-1}\circ ...\circ G_1]^{n_2}(\tilde{x}, \tilde{y})(v_x, v_y)\|\geq k \lambda^{n_2}
 $$
where  $k>0$ and $\lambda>1$. Since $ \mathcal{F}^j_{(\varepsilon_1, \varepsilon_2)}$ is a diffeomorphism for all $j=1, ..., qm+j$ and $[G_q \circ G_{q-1}\circ ...\circ G_1]^{n_2}(z)$ belongs to the compact set $\mathcal{A}'$, there exist constants $c,k>0$ and $\lambda>1$ such that:
\begin{eqnarray*}
 \|D \mathcal{F}_{(\varepsilon_1, \varepsilon_2)}^{qm+j} (z)v\| &=&\| D \mathcal{F}_{(\varepsilon_1, \varepsilon_2)}^{n_1} ([G_q \circ G_{q-1}\circ ...\circ G_1]^{n_2}(z))\|\times \|DG^{n_2}_{1, ..., q}(z)v\|
 \\
 &\geq& c \|D[G_q \circ G_{q-1}\circ ...\circ G_1]^{n_2}(z) v\|.\\
 &\geq& c k \lambda^{n_2}>0.
\end{eqnarray*}
 This finishes the proof of the lemma.
\end{proof}

  \section{Numerics}
  \label{s: numerics}
\begin{figure}[ht]
\begin{center}
\includegraphics[height=7.5cm]{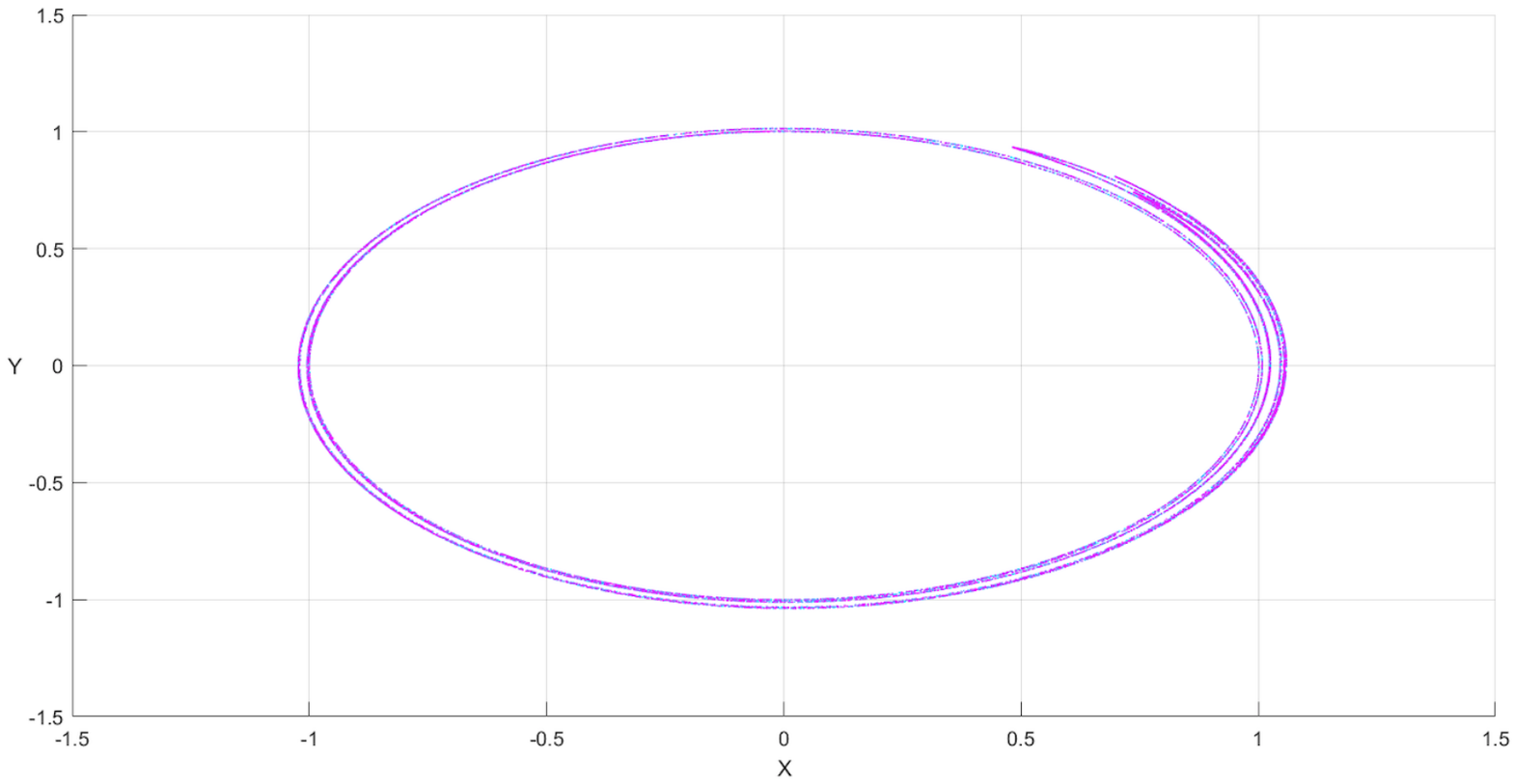}
\end{center}
\caption{\small Illustration of Theorem \ref{th:transv}. Orbit of $\mathcal{F}^n_{(\varepsilon_1, 0)}(x_0, y_0, t_0)$, for $1000<n<21000$, in rectangular coordinates $(X, Y)=(y\cos (x), y\sin (x))$.  Initial condition:  $(x_0, y_0, t_0)=(0.6961,1.3277, 0.5856)$, $\varepsilon_1 = 0.105$, $\delta_1 = 5$ and $\alpha_1= 6.2831= 2\alpha_2$. }
\label{Numerics A}
\end{figure}

\begin{figure}
\begin{center}
\includegraphics[height=7.5cm]{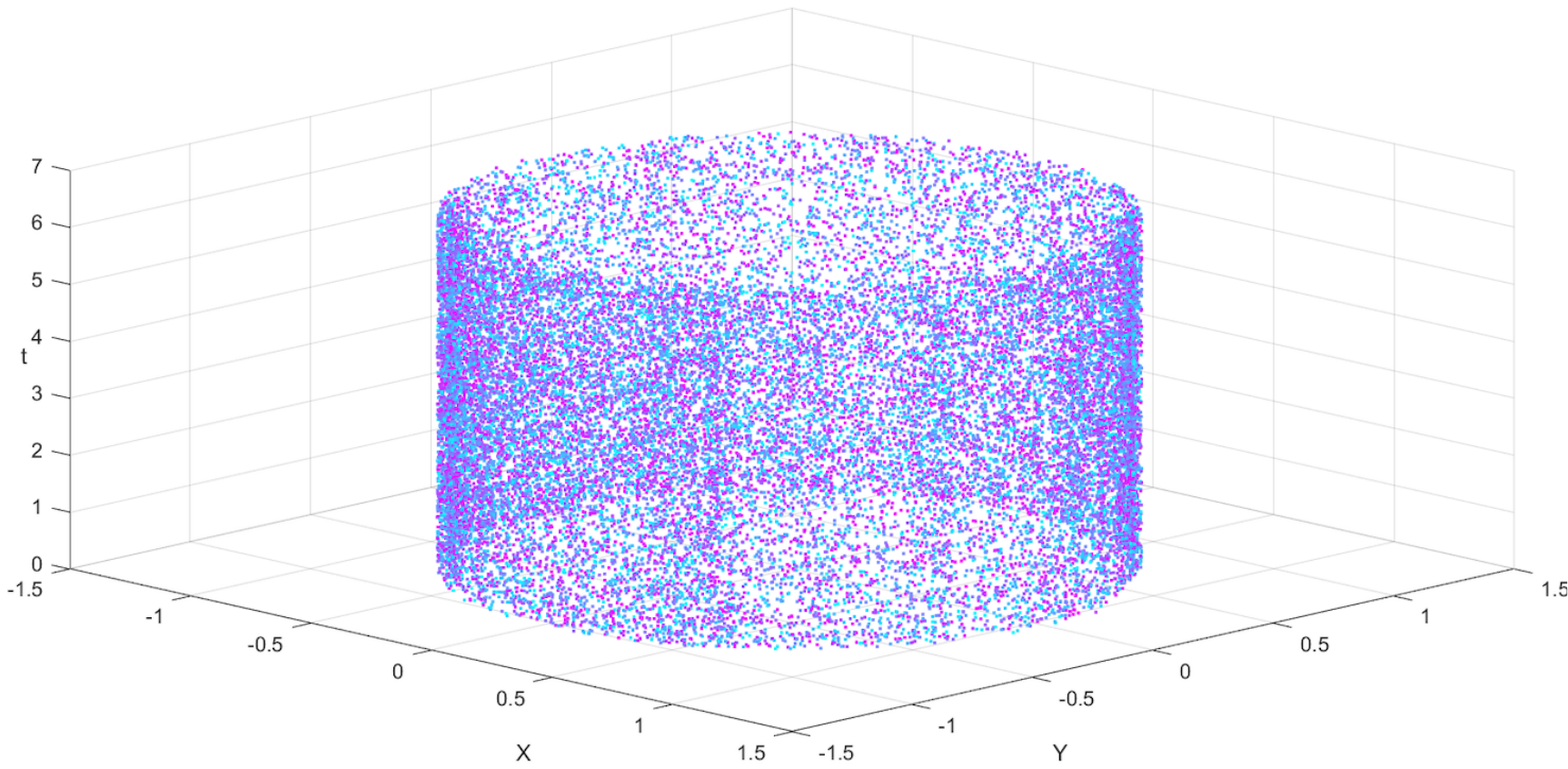}
\end{center}
\caption{\small Illustration of Theorem \ref{prop0}. Orbit of $\mathcal{F}^n_{(\varepsilon_1, 0)}(x_0, y_0, t_0)$, for $1500<n<31500$, in the coordinates $(X, Y, t)=(y\cos (x), y\sin (x), t)$. The planes defined by $t=0$ and  $t=2\pi$ are identified. Initial condition:  $(x_0, y_0, t_0)=(0.9073,1.4529, 0.5635)$, $\varepsilon_1 = 0.1$, $\delta_1 = 10$, $\delta_2= 0.001$, $\alpha_1= 6.2832$ and $\alpha_2= 4.4407$. }
\label{Numerics B}
\end{figure}

\begin{figure}
\begin{center}
\includegraphics[height=7.2cm]{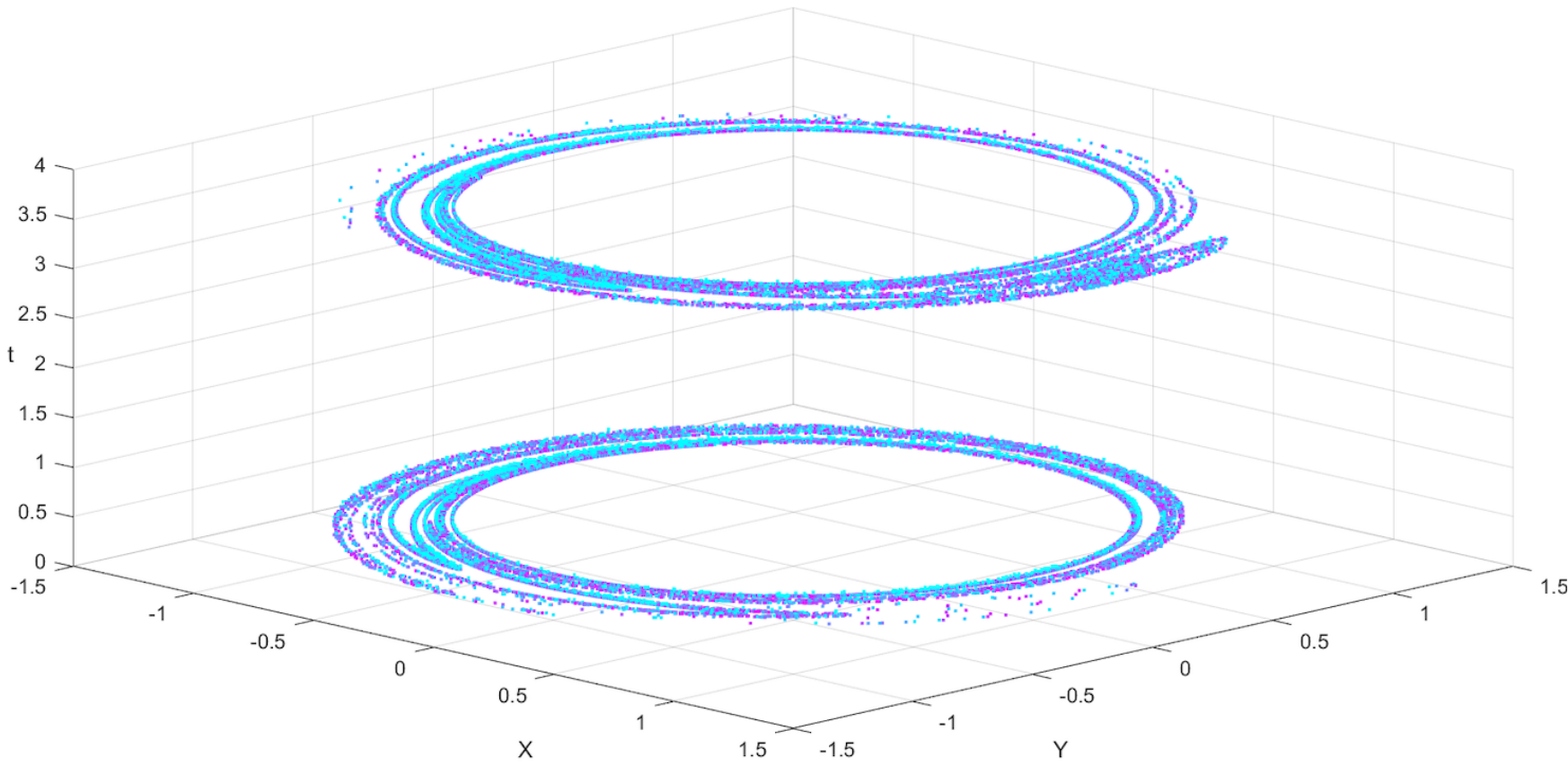}
\end{center}
\caption{\small Illustration of Theorem \ref{Th2.5}. Orbit of $\mathcal{F}^n_{(\varepsilon_1, \varepsilon_2)}(x_0, y_0, t_0)$, for $5000<n<105000$, in the coordinates $(X, Y, t)=(y\cos (x), y\sin (x), t)$.  Initial condition:  $(x_0, y_0, t_0)=(0.8394, 1.3789, 0.8716)$, $\varepsilon_1 = 0.2$,  $\varepsilon_2 = 0.1$, $\delta_1 = 5$, $\delta_2= 0.001$ and $\alpha_1= 6.2832=2\alpha_2$. Iterates of   $t_0$ under $\mathcal{F}_3$ are jumping along  a (stable) 2--periodic orbit.}

\label{Numerics D}
\end{figure}

\begin{figure}
\begin{center}
\includegraphics[height=7.2cm]{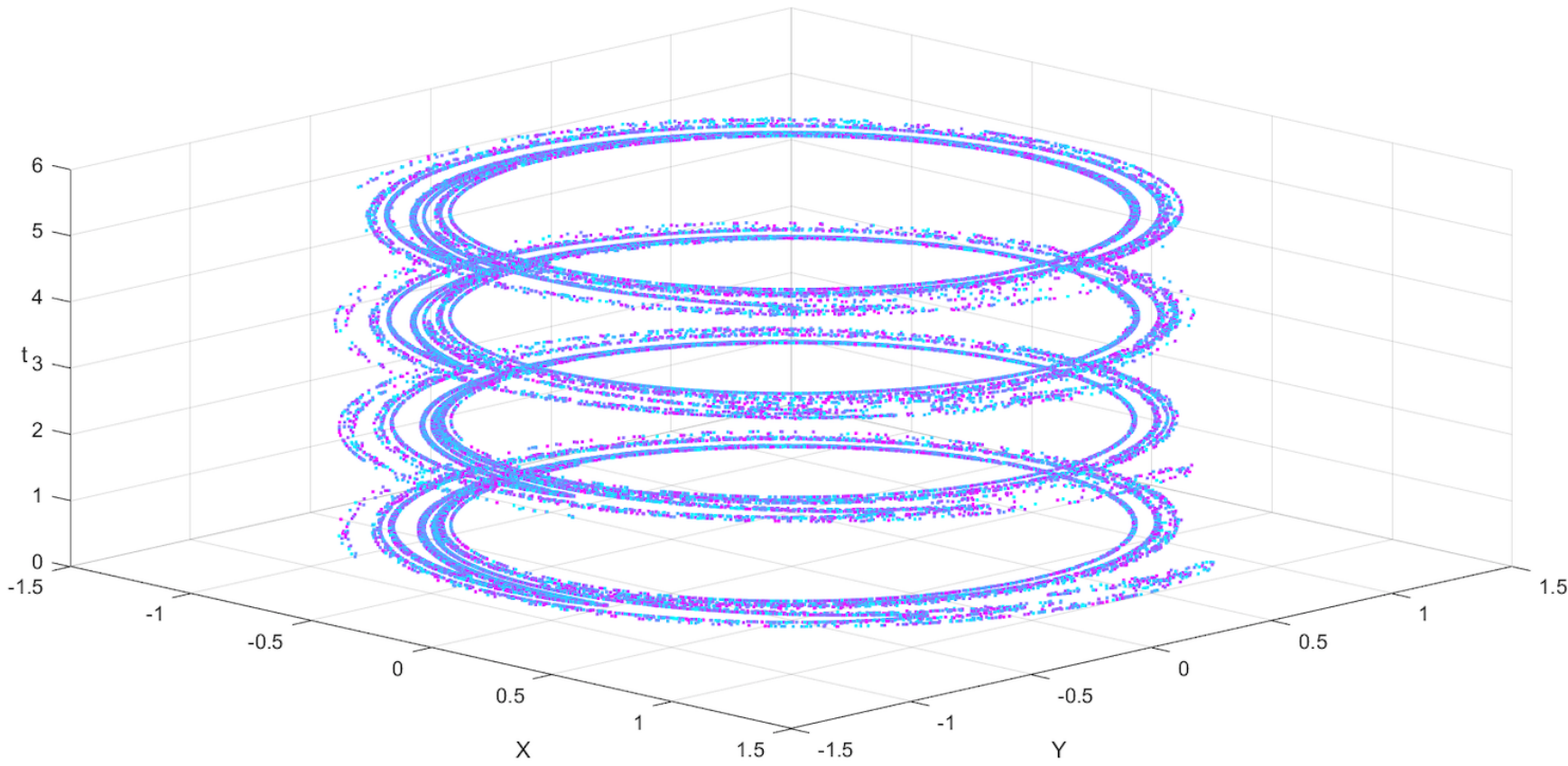}
\end{center}
\caption{\small Illustration of Theorem \ref{Th2.5}. Orbit of $\mathcal{F}^n_{(\varepsilon_1, \varepsilon_2)}(x_0, y_0, t_0)$, for $5000<n<105000$, in the coordinates $(X, Y, t)=(y\cos (x), y\sin (x), t)$.  Initial condition:  $(x_0, y_0, t_0)=(0.8162,1.0488, 0.6393)$, $\varepsilon_1 = 0.2$,  $\varepsilon_2 = 0.1$, $\delta_1 = 5$, $\delta_2= 1\times 10^{-7}$, $\alpha_1=6.2832$ and $\alpha_2=1.5708$. Iterates of   $t_0$ under $\mathcal{F}_3$ are jumping along  a (stable) 4--periodic orbit.}

\label{Numerics D2}
\end{figure}

  The aim of this section is to illustrate Theorems   \ref{th:transv}, \ref{prop0} and  \ref{Th2.5} through numerics using the software \emph{Matlab}(R2023b). All the simulations respect the dynamics of (\ref{map_general}).
Figures \ref{Numerics A}, \ref{Numerics B},  \ref{Numerics D} and  \ref{Numerics D2} depict a specific orbit   under the conditions $\delta=2$, $b=0.5$ and the following positive $C^3$--maps:
  \begin{eqnarray*}
 \Psi_1(x,y)= \Psi_2(x,y)=  g(x,y)= 1.1 + \sin(x)>0;\\
 \Psi_3(t)  = 1.1 + \sin(t)>0;\\
 \Psi_4 (x,y,t)= 1.1 + \sin(t)>0. 
\end{eqnarray*}

 These functions correspond to one of the simplest  positive $2\pi$--periodic maps respecting Hypotheses \textbf{(H3)--(H6)}.
  We use rectangular coordinates $X=y\cos x$ and $Y=y \sin x$ to allow further comparison with other works in the literature.
  
Roundoff errors are locally amplified by a large factor if $\delta_2$ is sufficiently wide. This can result in numerically observed behaviour which may be different from the actual dynamics, namely the existence of sinks.
The upper Lyapunov exponent (associated to the orbit under consideration) in all pictures are positive. 

 In Figure \ref{Numerics A}, we see the irreducible strange attractor $\Omega$ governing the  dynamics of $\mathcal{T}_{(\varepsilon_1, 0)}$.  It has one direction of instability and is structurally unstable, although their existence is a prevalent phenomenon   (a phenomenon already found in Figure 5 of \cite{BST98}). In Figure \ref{Numerics B}, we observe the existence of a quasi-periodic strange attractor in the uncoupled case $\mathcal{F}_{(\varepsilon_1, 0)}$. Note that the planes defined by $t=0$ and $t=2\pi$ are identified, so the image seems to cover a two--torus.
 In this case, the attractor may be seen as $\Omega \times \EU^1$. Finally, in Figures \ref{Numerics D} and  \ref{Numerics D2}, we  observe the dynamics of $\mathcal{F}_{(\varepsilon_1, \varepsilon_2)}$ where  $(\alpha_2, \delta_2)\in \mathcal{U}^{1/2}$ and $(\alpha_2, \delta_2)\in \mathcal{U}^{1/4}$, respectively. Although it seems to have several connected components, the  strange attractor is irreducible and displays a fractal structure. 

  \section{Discussion and concluding remarks}
  \label{s: discussion}
  This article studies the dynamics of the  two-parameter family of diffeomorphisms $\mathcal{F}_{(\varepsilon_1, \varepsilon_2)}$ (\ref{map_general}) defined on $\EU^1\times [1, 1+b]\times\EU^1$, with $b\gtrsim 0 $. This family  is obtained by the weak coupling of a rank-one family $\mathcal{T}_{(\varepsilon_1, 0)}$ and an Arnold diffeomorphism $\mathcal{F}_3$ on the circle $\EU^1$ , which is topologically conjugate to a rigid rotation. If $\delta_1\gg1$ and $\varepsilon_1\in \mathcal{G}_1$ (set of positive Lebesgue measure given by Theorem \ref{th:transv}), the  skew-product  \eqref{map_general}  has a component of chaotic behaviour and another with periodic or quasi-periodic motion.
 These families may behave periodically, quasi-periodically or chaotically, depending on specific character of the perturbation.

The main result   is the existence of rank-one strange attractors in the sense of \cite{WY2001, WY, WY2003, WY2008} for the family  $\mathcal{F}_{(\varepsilon_1, \varepsilon_2)}$. Making use of the existing theory, we have shown that the coupled diffeomorphism $\mathcal{F}_{(\varepsilon_1, \varepsilon_2)}$ has an invariant circle $\mathcal{C}$ of saddle type, such that  its stable and unstable manifolds   are bounded and   the orbits of all points in the absorbing domain of   $\mathcal{F}_{(\varepsilon_1, \varepsilon_2)}$ are attracted to the closure of its unstable manifold. The construction of $\mathcal{C}$ is described in Theorems \ref{prop0} and \ref{Th2.1}.

The main novelty used in this manuscript is the extension of the theory developed in \cite{WY}   to skew-products where one of the maps is Morse-Smale.
  In the case where the minimal set of $ \mathcal{F}_3$ is the entire circle, the coupling is governed by an irreducible quasi-periodic strange attractor  for a set of parameter values having positive   Lebesgue measure. Making use of Proposition \ref{Prop2.1WY}, Theorem \ref{Th2.5} couples  the dynamics  by using perturbations of  Misiurewicz-type maps.  This is valid because the ``simple'' dynamics of $\mathcal{F}_3$ is of the type Morse-Smale for $\delta_2$ satisfying $\textbf{(H5)}$.

 Similar results in the two-torus have been obtained by Jaeger \cite{Jager} where the author proved  the existence of quasi-periodic strange attractors/repellors in quasiperiodically forced circle maps under rather general conditions (stated in terms of $C^1$--estimates).  The strange attractors of \cite{Jager} carry the unique physical measure of the system, which determines the behaviour of Lebesgue-almost all initial conditions. 
 The results apply   to a forced version of the Arnold circle map.  
 
 
  Theorems \ref{prop0}, \ref{Th2.1} and  \ref{Th2.5}  will be  particularly useful  in the near future to study hyperchaos near a  network associated to a bifocus \cite{Barrientos_book, Rodrigues2020}.
 The full analysis of this scenario is not possible (by now) but the family $\mathcal{F}_{(\varepsilon_1, \varepsilon_2)}$ can be seen as a particular return map of the transition dynamics from chaos to persistent hyperchaos in \cite{Rodrigues2020}. 
 
Under the hypotheses of Theorem   \ref{Th2.5}, if $\delta_2>0$ satisfies \textbf{(H5)}, one observes either a quasi-periodic strange attractor or a strange attractor according to the position of  $(\alpha_2, \delta_2)$  in the resonant tongue. We conjecture that as $\delta_2$ increases, the deformation of the strange attractor shadowing $y=1$ is exaggerated, giving rise to rotational horseshoes created by double stretch-and-fold type actions.
 Iterates do not escape and wander around $\mathcal{A}'$. For $\delta_2 \gg 1$,  we guess that there exists a strange attractor with two positive Lyapunov exponents. The rigorous proof of this result is a challenge.

 The strange attractors in the present paper are nonuniformly hyperbolic and structurally unstable, although their existence is a persistent phenomenon.
 Techniques used  are valid for more general families   and also in the context of periodically perturbed systems. The study of the ergodic consequences of this dynamical scenario is the natural continuation of the present work.

\end{document}